%% file: equivFewOrbitToroids.tex
\documentclass[11pt]{article}
\usepackage[utf8]{inputenc}
\usepackage[english]{babel} %babel
\usepackage{hyperref}
\usepackage{amsmath, amsthm, amssymb, amsfonts} %ams packages
\usepackage[capitalise, noabbrev, nameinlink, poorman]{cleveref} %nice references
\usepackage{graphicx, color}%paquetes usados en caso de usar gráficos
\usepackage[figurewithin=none]{caption,subcaption} %option for change numering of figures
\usepackage{enumitem} %personalized lists 
\usepackage{marginnote} %personalized marginotes
\usepackage[]{geometry} %to change margins
\usepackage{soul} % For using "\st" command
\usepackage{array} %bmatrix and similar stuff
\usepackage{filecontents} %to use embedded bibtex

\usepackage{tikz-cd}
\usepackage{todonotes}
\usepackage{multirow} 

\renewcommand{\leq}{\leqslant} %menores y mayores iguales decentes
\renewcommand{\geq}{\geqslant}
\renewcommand{\epsilon}{\varepsilon} %epsilon decente
\renewcommand{\subset}{\subseteq} %subseteq decente
\renewcommand{\{}{\lbrace}
\renewcommand{\}}{\rbrace}
\newcommand{\sm}{\setminus} %setminus corto

              \def\cO{{\mathcal{O}}} \def\cP{{\mathcal{P}}}      \def\cU{{\mathcal{U}}}    

    \def\bE{{\mathbb{E}}}         \def\bN{{\mathbb{N}}}            \def\bZ{{\mathbb{Z}}}

\theoremstyle{plain}
\newtheorem{theorem}{Theorem}
\newtheorem{lemma}[theorem]{Lemma}

\newtheorem{corollary}[theorem]{Corollary}

\theoremstyle{definition}

\theoremstyle{remark}

\newcommand{\E}[1][n]{\bE^{#1}}

\newcommand{\LL}{\Lambda}
\newcommand{\bLL}{\mathbf{\LL}}

\newcommand{\cln}[1][n-1]{\Lambda_{(1, 0^{#1})}}
\newcommand{\fcln}[1][n-2]{\Lambda_{(1,1, 0^{#1})}}
\newcommand{\bcln}[1][n]{\Lambda_{(1^{#1})}}

\newcommand{\clFour}{\Lambda_{(1, 0,0,0)}}
\newcommand{\fclFour}{\Lambda_{(1,1, 0,0)}}
\newcommand{\bclFour}{\Lambda_{(1,1,1,1)}}

\newcommand{\latEven}{\Lambda_{0}}
\newcommand{\latOdd}{\Lambda_{1}}
\newcommand{\regLat}[1][n]{\Lambda_{(1^{k},0^{#1-k})}}
\newcommand{\LatProd}{\Lambda_{(1,1)} \times \Lambda_{(1,1)}}

\newcommand{\trasLat}[1][w]{\bigcup_{k \in \bZ} \left(\LL_{0} + k#1 \right)}

\def\minid{\chi}
\def\triv{\mathbf{1}}
\def\id{id}

\newcommand{\symn}[1][n]{\mathbf{S_{#1}}}
\newcommand{\altn}[1][n]{\mathbf{A_{#1}}}
\newcommand{\abeln}[1][n]{\mathbf{C^{#1}_{2}}}
\newcommand{\abelnplus}[1][n]{\mathbf{(C^{#1}_{2})^{+}}}
\newcommand{\died}{\mathbf{D_{4}}}

\newcommand{\stabcube}[1][n]{\abeln[#1] \rtimes \symn[#1]}
\newcommand{\pgl}{\mathbf{PGL_{2}(5)}}

\newcommand{\UoverL}{\cU/\bLL}

\DeclareMathOperator{\G}{G}
\DeclareMathOperator{\T}{T}
\DeclareMathOperator{\isom}{Isom}
\DeclareMathOperator{\Go}{\G_{o}}
\DeclareMathOperator{\normI}{Norm_{\isom(\E)}}
\DeclareMathOperator{\normGU}{Norm_{G(\cU)}}
\DeclareMathOperator{\aut}{Aut}

\newcommand{\autUoverL}{\aut(\cU/\bLL)}

\begin{document}

\title{Equivelar toroids with few flag-orbits \thanks{Supported by PAPIIT UNAM project grant IN101615.}}

\author{José Collins \\ collins@ciencias.unam.mx \\ Instituto de Matemáticas UNAM \\ Mexico City, México. \and Antonio Montero \\ amontero@matmor.unam.mx\\  Centro de Ciencias Matemáticas UNAM \\ Morelia, México. } %thanks tero: CONACyT grant 414098 
%\date{}
\maketitle

\begin{abstract}
	\input{abstract}
\end{abstract}

%  \listoftodos
\section{Introduction} \label{sec:intro}
\input{intro}

\section{Basic notions} \label{sec:basics}
\input{toroids}
%yo creo que es buena idea meter también unas tablitas con los generadores de las teselaciones, muy como las tienen estos cuates para las de dimensión 4
% muy parecido a HPOW, incluir resultados de esos weyes y de Egon, Peter y Michael Hartley.

\section{Few-orbit cubic toroids} \label{sec:cubic}
	\input{cubic}

	%Un poquito de estructura del grupo
	% normales del abeliano 
	% subgrupos de indice 2, todos continen al (C_{2}^{n})^{+} \rtimes A_n
	% Latices invariantes bajo (C_{2}^{n})^{+} \rtimes A_n
	% clasificación de toros de 2 órbitas
	% el único grupo de indice k con 2 < k < n es (C_{2}^{n})^{+} \rtimes A_n
	% clasificación de toros de k órbitas con 2 < k < n 
	% el unico grupo de indice n en S_{n} es S_{n-1} <- Tero cree que ya valió pito y debe estudiar (Dixon)
	% Clasificación de toros de n órbitas. <- Salvo la valida de pito de arriba, ya están. 
\section{Few-orbit toroids of type \texorpdfstring{$\{3,3,4,3\}$}{3,3,4,3} and \texorpdfstring{$\{3,4,3,3\}$}{3,4,3,3}} \label{sec:nocubic} 
	% Un poquito de estructura del grupo
	% clasificación de clases de conjugación de índice 2,3,4 <- No ha salido. Hasta ahora sólo podemos apelar a GAP
	% Clasificación de few-orbits toroids.
	\input{nonCubic}
	\section{Symmetry type of few-orbit toroids} \label{sec:symtype} 
	\input{symType}

% \section{Conclusions}\label{sec:conc}
	%Hacer teoremon.
%Content

\bibliographystyle{plain}
\bibliography{equivFewOrbitToroids.bib}

\end{document}

%% file: abstract.tex
	An $(n+1)$-toroid is a quotient of a tessellation of the $n$-dimensional Euclidean space with a lattice group. Toroids are generalizations of maps in the torus on higher dimensions and also provide examples of abstract polytopes. Equivelar toroids are those that are induced by regular tessellations. In this paper we present a classification of equivelar $(n+1)$-toroids with at most $n$ flag-orbits; in particular, we discuss a classification of $2$-orbit toroids of arbitrary dimension.

%% file: intro.tex
The study of symmetric discrete objects on both combinatorial and geometrical approach has been of interest in recent years. In particular, symmetric maps on surfaces has been one of the most studied topics.

In \cite{coxeterMoserGenandRelforDG} Coxeter and Moser present the classification of regular (reflexible) and chiral (regular irreflexible) maps on the torus. All such maps arise as quotients of a regular tessellation of the euclidean plane. Several results regarding a classification of highly symmetric maps on surfaces of small genus have been obtained since. See \cite{conderDobcsanyi_RegMapsSmallGenus,conderRegMapsChar1to200}.

When looking for generalizations of maps to higher dimensions the approach of abstract polytopes has been one of the most studied. The number of recent contributions to the theory of symmetric abstract polytopes is large. Many of such results may be found in \cite{ARP}.

One concept that generalizes maps on a combinatorial way while keeping the topological idea behind maps is that of tessellations of space forms (see \cite[Chapter 6]{ARP}). Euclidean space forms are probably the most studied on the setting of symmetric tessellations. Among those, the $n$-dimensional torus is probably the most well understood.

When talking about symmetric structures on the $n$-dimensional torus, much of the work follow the ideas introduced by Coxeter and Moser in \cite{coxeterMoserGenandRelforDG}. \emph{Toroids} are generalizations to higher dimension of maps on the torus and may be regarded as tessellations of the $n$-dimensional torus. 

Several classification results of highly symmetric toroids have been developed in the recent years. In \cite{mcMullenShulte_HigherToroidalRP} McMullen and Schulte classify regular toroids of arbitrary dimension. They also show that there are no chiral toroids of dimension higher than $3$. In \cite{harleyMcMullenSchulte_SymTessESF} Hartley, McMullen and Schulte extend this result and show that the only Euclidean space form that admits chiral tessellations is the $n$-torus and prove that this is only possible when $n=2$. 

In \cite{brehmKuhnel_EquivelarMapsTorus} Brehm and Kühnel classify the equivelar maps in the two dimensional torus. This classification was also achieved by Hubard et al. in \cite{hubardOrbanicPellicerWeissEquiv4Toroids}. They also extend the techniques to classify equivelar tessellations of the $3$ dimensional torus (rank-$4$ toroids). 

As a consequence of their results in \cite{hubardOrbanicPellicerWeissEquiv4Toroids}, Hubard, Orbanić, Pellicer and Weiss found that there are no equivelar $2$-orbit $(3+1)$-toroids. Therefore, they mentioned that even though a classification of $(n+1)$-toroids for $n > 4$ seems to be too hard to achieve with their techniques, it would be of interest to obtain a classification of $2$-orbit toroids of rank $n+1$ for $n \geq 4$.

In this paper we give a classification of equivelar $(n+1)$-toroids with at most $n$ orbits. In particular, $2$-orbit $(n+1)$-toroids are classified for arbitrary dimension. If $n\leq 3$ the classification is a consequence of the results in \cite{coxeterMoserGenandRelforDG} and \cite{hubardOrbanicPellicerWeissEquiv4Toroids}. The main results of this article can be summarized in the following theorem.

\begin{theorem} \label{thm:toroids}
	Let $n\geq 4$. The classification of equivelar $(n+1)$-toroids is explained in Table\nobreakspace \ref {tab:classification}.
\end{theorem}

\begin{table}
 \begin{center}
 \resizebox{\textwidth}{!}{%
\newcommand{\mc}[3]{\multicolumn{#1}{#2}{#3}}
\begin{tabular}{|c|m{0.38\textwidth}|m{.38\textwidth}|}
 \cline{2-3}
\mc{1}{c|}{}& \mc{1}{c|}{$n=4$} & \mc{1}{c|}{$n\geq4$}\\ \hline
 \multirow{8}{0.15\textwidth}{$\{4,3^{n-2},4\}$} & Three families of regular toroids. & Three families of regular toroids.  \\ \cline{2-3}
& One family of $2$-orbit toroids in class $2_{\{1,2,3\}}$& If $n$ is odd, there are no $2$-orbit toroids. If $n$ is even there is one family of torois in class $2_{\{1,2, \dots, n-1\}}$.\\ \cline{2-3}
& One family of $3$-orbit toroids. & There are no $k$-orbit toroids if $2 < k < n$. \\ \cline{2-3}
& Five infinite families of toroids with $4$ flag-orbits, all with the same symetry type. & Five infinite families of toroids with $n$ flag-orbits, all with the same symetry type. \\ \hline
	\multirow{6}{0.15\textwidth}{$\{3,3,4,3\}$ and $\{3,4,3,3\}$} 	& Two families of regular toroids. & \mc{1}{c}{} \\ \cline{2-2}
& One family of $2$-orbit toroids in class $2_{\{3,4\}}$ (or dually, in class $2_{\{0,1\}}$). & \mc{1}{c}{}\\ \cline{2-2}
& Two families of $3$-orbit toroids with different symmetry type. & \mc{1}{c}{}\\ \cline{2-2} 
& There are no $4$-orbit toroids.& \mc{1}{c}{}\\ \cline{1-2}
\end{tabular}
 }
\caption{Classification of equivelar few-orbit toroids.} \label{tab:classification}
\end{center}
\end{table}

% For $n \geq 4$ the classification of equivelar $(n+1)$-toroids with at most $n$ flag-orbits is described in the following list.
% 	\begin{itemize}
% 		\item For every $n \geq 4$ there are three infinite families of regular toroids of type $\{4, 3^{n-2}, 4\}$.
% 		\item Equivelar $2$-orbit toroids of rank $(n+1)$ and type $\{4, 3^{n-2},4\}$ exist if and only if $n$ is even. Such toroids belong to class $2_{\{1, \dots, n-1\}}$.
% 		\item For $2 < k < n$, there are no equivelar toroids  of type $\{4, 3^{n-2}, 4\}$ with $k$ flag-orbits .
% 		\item There are five infinite families of $(n+1)$-toroids of type $\{4, 3^{n-2}, 4\}$ whose automorphism group has $n$ flag-orbits. All the families share the symmetry type.
% 		\item If $n=4$, there are two  families of regular toroids of type $\{3,3,4,3 \}$ and two families of regular toroids of type $\{3,4,3,3\}$.
% 		\item There is one family of $2$-orbit toroids of type $\{3,3,4,3 \}$ and one family of $2$-orbit toroids of type $\{3,4,3,3\}$. Those families belong to classes $2_{\{3,4\}}$ and $2_{\{0,1\}} $, respectively.
% 		\item There exist two families of equivelar toroids of type $\{3,3,4,3 \}$ and two families of toroids of type $\{3,4,3,3\}$ with $3$ flag-orbits. Those families have different symmetry type.
% 		
% 	\end{itemize}

%% file: toroids.tex
%yo creo que es buena idea meter también unas tablitas con los generadores de las teselaciones,
%muy como las tienen estos cuates para las de dimensión 4
% muy parecido a HPOW, incluir resultados de esos weyes y de Egon, Peter y Michael Hartley.

%Tessellations
%symmetries of tesselations
%regular tessellations

\subsection{Tessellations of the Euclidean Space} \label{sec:tessellations}
In this section we introduce basic concepts about Euclidean tessellations and toroids. We focus mainly on those with high degree of symmetry. Readers interested in further details are referred to \cite[Chapter 6]{ARP} and \cite{hubardOrbanicPellicerWeissEquiv4Toroids}.

A \emph{convex $n$-polytope} $\cP$ is the convex hull of a finite set of points of $\E$ such that the interior of $\cP$ is non-empty. If $\cP$ is a convex $n$-polytope, then $F \subset \cP$ is a \emph{face} of $\cP$ if $F = \cP \cap \Pi$ for some hyperplane $\Pi$ that leaves $\cP$ contained in one of the closed half spaces determined by $\Pi$. If the affine dimension of a face $F$ is $i$ for some $i \in \{0, \dots , n-1\}$, then $F$ is an \emph{$i$-face}. The $0$-faces, $1$-faces and $(n-1)$-faces of a convex $n$-polytope are also called \emph{vertices}, \emph{edges} and \emph{facets}, respectively. We usually consider a convex $n$-polytope $\cP$ itself as its (unique) $n$-face. Readers interested in more details concerning basic notions of convex polytopes are referred to \cite[Chapter 5]{matousek_DiscreteGeometry}.

A \emph{tessellation} of the Euclidean space $\E$ (or an \emph{Euclidean tessellation}) is a family $\cU$ of convex $n$-polytopes that is \emph{locally finite}, meaning that every compact set of $\E$ meets only finitely many members of $\cU$. We also require that $\cU$ covers the space and tile it in a face-to-face manner. This is, every two members of $\cU$ that have non-empty intersection have disjoint interiors, and they meet in a common $i$-face for some $ i \in \{0, \ldots, n-1\}$. If $\cU$ is a Euclidean tessellation, the elements of $\cU$ are called \emph{cells}.

 A flag in a convex $n$-polytope is an $(n+1)$-tuple of incident faces containing exactly one face of each dimension, including the polytope itself. This definition extends naturally to tessellations of $\E$. A \emph{flag} of a tessellation $\cU$ is a flag of any cell of $\cU$. It is sometimes useful to identify a flag $\Phi$ of $\cU$ with the non regular $n$-simplex induced by the centroids of the faces of $\Phi$. Observe that given a flag $\Phi$ of $\cU$ and $i \in \{0, \dots, n\}$, there exists exactly one flag $\Phi^{i}$ of $\cU$ that differs from $\Phi$ only in the face of dimension $i$. In this situation we say that $\Phi$ and $\Phi^{i}$ are \emph{adjacent} (or \emph{$i$-adjacent} if we want to emphasize on $i$). 
 
 A \emph{symmetry} of a tessellation $\cU$ is an isometry of $\E$ that preserves $\cU$. The \emph{group of symmetries of $\cU$} is denoted by $\G(\cU)$. It is not hard to see that $\G(\cU)$ acts freely on the set of flags of $\cU$. A tessellation $\cU$ is \emph{regular} if the action of $\G(\cU)$ on the flags of $\cU$ is transitive.
 
 The \emph{dual tessellation} of a regular Euclidean tessellation $\cU$, usually denoted by $\cU^{\ast}$, is the tessellation whose cells are the polytopes given by the convex hull of the centroids of the cells of $\cU$ incident to a common vertex of $\cU$. A tessellation is \emph{self-dual} if it is isometric to its dual.

The \emph{Schläfli type} (or \emph{type}, for short) of a convex $n$-polytope $\cP$ is defined recursively as follows. If $n=2$ then $\cP$ is a convex $p$-gon for some $p$ and we say that $\cP$ has Schläfli type $\{p\}$. For $n \geq 3$, whenever all the facets of a convex $n$-polytope have type $\{p_{1}, \dots, p_{n-2}\}$ and there are exactly $p_{n-1}$ facets around each $(n-3)$-face of $\cP$, we say that $\cP$ has Schläfli type $\{p_{1}, \dots, p_{n-1}\}$. Observe that not every convex polytope has a well-defined Schläfli type, however all regular convex polytopes do. The notion of Schläfli type extends to tessellations in a natural way. We say that $\cU$ has (Schläfli) type $\{p_{1}, \dots p_{n}\}$ if all the cells of $\cU$ have type $\{ p_{1}, \dots, p_{n-1} \}$ and the number of cells around each $(n-2)$-face of $\cU$ is $p_{n}$. 
% 
% In this sense, the cube has Schläfli type $\{4,3\}$ and the icosahedron has Schläfli type $\{3,5\}$.

Regular tessellations are well known. There exists a self-dual regular tessellation with cubes on $\E$ with type $\{4, 3^{n-2} ,4\}$. Here $3^{n-2}$ denotes a sequence of $3$ with length $n-2$; if there is no possible confusion, in this work we will use exponents to denote a sequence of equal symbols. In $\E[2]$ there exists a regular tessellation with equilateral triangles and type $\{3,6\}$, and a tessellation with regular hexagons and type $\{6,3\}$. Those two are dual of each other. In $\E[4]$ there is another pair of regular tessellations, one with 24-cells as facets and type $\{3,4,3,3\}$ and its dual of type $\{3,3,4,3\}$ whose cells are four dimensional cross-polytopes. These tessellations are unique up to similarity and they complete the list of regular tessellations of the Euclidean $n$-space \cite[Table II]{coxeterRegularPolytopes}. In Table\nobreakspace \ref {tab:tessellations} we give explicit coordinates for the vertex set of one of each pair of dual regular tessellations. Those coordinates determine uniquely the tessellation.

%distinguished generators and coxeter groups.

If $\cU$ is a regular tessellation of $\E$ and $\Phi$ is a fixed \emph{base flag}, there exist $R_{0}, \dots , R_{n}$  symmetries of $\cU$ such that $\Phi R_{i} = \Phi^{i}$. The symmetries $R_{0}, \dots , R_{n}$ are the reflections on the facets of the $n$-simplex induced by $\Phi$ and generate the group $\G(\cU)$. If $\cU$ has type $\{p_{1}, \dots, p_{n}\}$, the group $\G(\cU)$ with the generators $R_{0}, \dots , R_{n}$ is the string Coxeter group $[p_{1}, \dots, p_{n}]$ meaning that \[\G(\cU) = \left\langle R_{0}, \dots , R_{n} \ |\  (R_{i}R_{j})^{p_{i,j}}=id\right\rangle\] is a presentation for the group, where $p_{i,j} = p_{j,i}$, $p_{i,i}=2$, $p_{i,j}=2$ if $|i-j| > 1$, and $p_{i-1,i}=p_{i}$. Moreover, the group of symmetries of the cell contained in the base flag is $\langle R_{0}, \dots , R_{n-1} \rangle \cong [p_{1}, \dots, p_{n-1}]$ and the stabilizer of the vertex of $\Phi$ is  $\langle R_{1}, \dots , R_{n} \rangle \cong [p_{2}, \dots, p_{n}]$ (see \cite[section 3A]{ARP} for details). 

Note that if $R_{0}, \dots, R_{n}$ denote the distinguished generators of $\cU$ with respect to a base flag $\Phi$, then $R_{n}, \dots, R_{0}$ act as distinguished generators of $\cU^{\ast}$ with respect to some flag. A consequence of this is that $\G(\cU) = \G(\cU^{\ast})$. In Table\nobreakspace \ref {tab:tessellations} we give coordinates and explicit expressions for $R_{0}, \dots, R_{n}$ for one of each dual pair of regular tessellations.

\begin{table}
\begin{center}

\resizebox{\textwidth}{!}{%
\newcommand{\mc}[3]{\multicolumn{#1}{#2}{#3}}
\begin{tabular}{|c|c|l|} 
\hline
\mc{1}{|c}{\multirow{2}{*}{$\cU$}} & \mc{1}{|c|}{\multirow{2}{*}{Vertex set of $\cU$}} & \mc{1}{c|}{Generators of $\G(\cU)$}\\ 
& & \mc{1}{c|}{$R_{i}: (x_{1}, \dots x_{n}) \mapsto $}\\ \hline
\multirow{3}{*}{$\{3,6\}$} & \multirow{3}{*}{$\left\{a(1,0) + b \left(\frac{1}{2}, \frac{\sqrt{3}}{2}\right) \ : \ a,b \in \bZ \right\}$} & $(1-x_{1},x_{2})$ if $i=0$, \\
 &  & $\left(\frac{1}{2} x_{1}+ \frac{\sqrt{3}}{2}x_{2}, \frac{\sqrt{3}}{2}x_{1} -  \frac{1}{2}x_{2} \right)$ if $i=1$,\\
 &  & $(x_{1},-x_{2})$ if $i=2$.\\ \hline
\multirow{3}{*}{$\{4,3^{n-2},4\}$} &  \multirow{3}{*}{$\bZ^{n}$} & $(1-x_{1},x_{2}, \dots ,x_{n})$ if $i=0$, \\
 &  & $(x_{1}, \dots, x_{i-1}, x_{i+1}, x_{i}, \dots , x_{n})$ if $1 \leq i < n$, \\
 &  & $(x_{1}, \dots, -x_{n})$ if $i=n$. \\ \hline
 \multirow{6}{*}{$\{3,3,4,3\}$} &  \multirow{6}{*}{$\left\{(x_{1},x_{2}, x_{3},x_{4}) \in \bZ^{4} \ : x_{1}\equiv x_{2} \equiv x_{3} \equiv x_{4} \pmod{2}\ \right\}$} & $(2-x_{1},x_{2}, \dots ,x_{n})$ if $i=0$, \\
 &  & $(x,x- x_{3}- x_{4}, x-x_{2}-x_{4},x-x_{2}- x_{3})$ \\
 &  & where $x=\frac{1}{2}(x_{1} + x_{2}+x_{3}+x_{4})$, if $i=1$, \\
 &  & $(x_{1}, x_{2}, x_{3}, - x_{4})$ if $i=2$, \\ 
 &  & $(x_{1}, x_{2}, x_{4},  x_{3})$ if $i=3$, \\ 
 &  & $(x_{1}, x_{3}, x_{2},  x_{4})$ if $i=4$. \\ 
\hline
\end{tabular}
}%
\caption{Vertex-set and distinguished generators of regular tessellations.} \label{tab:tessellations}
\end{center}
\end{table}

\subsection{Equivelar toroids} \label{sec:toroids}

Toroids generalize maps in the $2$-dimensional torus to higher dimensions. In this section we will discuss the basic results about toroids and their symmetries, with special interest in those induced by regular tessellations of $\E$.

Let $0 \leq d \leq n$, a \emph{rank-$d$ lattice group} in $\E$ is a group generated by $d$ translations with linearly independent translation vectors. If $\bLL = \langle t_{1} , \dots , t_{d} \rangle$ is a  lattice group and $v_{i}$ is the translation vector of $t_{i}$, then the \emph{lattice} $\LL$ induced by $\bLL$ is the orbit of the origin $o$ under $\bLL$, that is \[\LL = o \bLL = \{a_{1}v_{1} + \cdots + a_{d}v_{d} \ \colon \ a_{1}, \dots , a_{d} \in \bZ\}.\] In this case we say that $\{v_{1}, \dots, v_{d}\}$ is a \emph{basis} for $\LL$.

We denote by $\cln$ the \emph{cubic lattice}, which consists of all the points of $\E$ with integer coordinates with respect to the standard basis $\{e_{1}, \dots e_{n}\}$ of $\E$. Observe that $\{e_{1}, \dots e_{n}\}$ is also a basis for $\cln$. The lattice $\fcln$ is the rank-$n$ lattice consisting of the points of integer coordinates of $\E$ whose coordinate sum is even. A basis for $\fcln$ is given by $\{2e_{1}, e_{2}-e_{1}, e_{3}-e_{2}, \dots, e_{n}-e_{n-1}\}$. We denote by $\bcln$ the lattice consisting of the points whose coordinates are all integers having the same parity. A basis for this lattice is $\{2e_{1}, \dots, 2e_{n-1}, e_{1}+ \cdots + e_{n}\}$. Note that $\cln$, $\fcln$ and $\bcln$ are contained in the vertex set of $\{4, 3^{n-2}, 4\}$ and $\bclFour$ is precisely the vertex set of $\{3,3,4,3\}$. The corresponding lattice groups are denoted by $\mathbf{\cln}$, $\mathbf{\fcln}$ and $\mathbf{\bcln}$. Note that this notation is slightly different from that of \cite{ARP}. 

Let $\T(\cU)$ denote the group of translations of an Euclidean tessellation $\cU$. A \emph{toroid of rank $(n+1)$} or \emph{$(n+1)$-toroid} is the quotient of a tessellation $\cU$ of $\E$ by a rank-$n$ lattice group $\bLL \leq \T(\cU)$. We say that $\bLL$ \emph{induces} the toroid, and denote the latter by $\cU/\bLL$. If $\cU$ is regular of type $\{p_{1}, \dots , p_{n}\}$ we say that the toroid is \emph{equivelar} of (Schläfli) type $\{p_{1}, \dots , p_{n}\}$; in this situation we also denote the toroid induced by $\bLL$ by $\{p_{1}, \dots , p_{n}\}_{\bLL}$ (cf. \cite[Chapter 7]{coxeterMoserGenandRelforDG} and \cite[Chapter 6]{ARP}). An $(n+1)$-toroid may be regarded as a tessellation of the $n$-dimensional torus $\E / \bLL$.

For all regular tessellations $\cU$ except $\{6,3\}$ and $\{3,4,3,3\}$ the group of translations $\T(\cU)$ acts transitively on the vertex set of $\cU$. Therefore the vertex set of $\cU$ may be identified with the lattice associated to $\T(\cU)$ and the group of automorphisms $\G(\cU)$ is of the form $\T(\cU)\rtimes \Go(\cU)$, where $\Go(\cU)$ denotes the stabilizer of the origin $o$ (see for example \cite[Chapter 6]{ARP}). For now on, we restrict our study to regular tessellations whose vertex-set is a lattice. The results regarding toroids of type $\{6,3\}$ and $\{3,4,3,3\}$ may be recovered by duality.

Let $\isom(\E)$ denote the group of isometries of $\E$. If $t_{v}$ is the translation by a vector $v$ and $S \in \isom(\E)$  fixing the origin $o$, then $S^{-1} t_{v} S = t_{vS}$, the translation by $vS$. In other words, if $\LL$ is the lattice associated to $\bLL$, then $\LL S$ is the lattice associated to $S^{-1} \bLL S$. Therefore if there exists an isometry mapping a lattice $\LL$ to another lattice $\LL'$, then there exists an isometry $S$ that fixes $o$ that maps $\LL$ to $\LL'$. In this case the corresponding tori $\E / \bLL$ and $\E /\bLL'$ are isometric. Geometrically this means that $S$ maps fundamental regions of $\bLL$ to fundamental regions of $\bLL'$. This implies that two toroids $\cU/\bLL$ and $\cU/\bLL'$ are isometric if and only if there exists an isometry $S \in \Go(\cU)$ mapping $\LL$ to $\LL'$ or, equivalently $S^{-1} \bLL S = \bLL'$. 

With the notation given above, when $\bLL = \bLL'$, an isometry $S$ of $\E$ induces an isometry $\overline{S}$ that makes the diagram \eqref{eq:diagtorus} commutative if and only if $S$ normalizes $\bLL$. Furthermore, two isometries of $\E$ induce the same isometry of $\E / \cU$ if and only if they  differ by an element of $\bLL$. In particular, all the elements of $\bLL$ induce a trivial isometry of $\E/\bLL$. This implies that the group $\normI(\bLL)/\bLL$ acts as a group of isometries of $\E/\bLL$. It can be proven that every isometry of $\E/\bLL$ is given this way, that is $\isom(\E/\bLL) \cong \normI(\bLL)/\bLL$ (see \cite[p.336]{ratcliffeHypManifolds} and \cite[Section 6A]{ARP}).

\begin{equation} \label{eq:diagtorus}
\begin{tikzcd}
\E \arrow{r}{S} \arrow{d}{} & \E \arrow{d}{} \\
\E/\bLL \arrow[dashed]{r}{\overline{S}} & \E / \bLL
\end{tikzcd}
\end{equation}

With the previous discussion in mind, it makes sense to define the group of automorphisms of a toroid $\UoverL$, denoted by $\autUoverL$, as the quotient $\normGU(\bLL)/\bLL$. Intuitively speaking, $\normGU(\bLL)$ denotes the symmetries of $\cU$ that are compatible with the quotient by $\bLL$. In the same sense, two toroids $\UoverL$ and $\UoverL'$ are \emph{isomorphic} if $\bLL$ and $\bLL'$ are conjugate in $\isom(\E)$. 

If an isometry $S$ normalizes $\bLL$ then we say that $S$ \emph{induces} or \emph{projects to} an automorphism of $\UoverL$ (namely, to the automorphism $\bLL S \in \normGU(\bLL)/\bLL$). If $S' \in \Go(\cU)$, then $S'$ normalizes $\bLL$ if and only if $S'$ preserves $\LL$. Since every element $ S \in \G(\cU)$ may be written as a product $tS'$ with $t \in \T(\cU)$ and $S' \in \Go(\cU)$ and every translation normalizes $\bLL$, an isometry $S$ induces an automorphism of $\cU/\bLL$ if and only if $S'$ preserves $\LL$. Therefore we may restrict our analysis to elements of $\Go(\cU)$.

The $i$-faces of a toroid $\cU/\bLL$ are the orbits of the $i$-faces of $\cU$ under $\bLL$. Whenever all the vertices on each cell of $\cU$ are different under the action of $\bLL$ the set of faces of $\cU/\bLL$ has the structure of an abstract polytope (in the sense of \cite{ARP}). In this case, the symmetry properties of $\UoverL$ as an abstract polytope  coincide with those as a toroid.  However, even when $\UoverL$ is not an abstract polytope, we may define the set of flags of $\cU/\bLL$ as the set of orbits of flags of $\cU$ under $\bLL$. The group $\aut(\cU/\bLL)$ has a well defined action on the set of flags of $\cU/\bLL$ by $(\Phi\bLL) \bLL S = \Phi S \bLL$ with $\Phi$ a flag of $\cU$ and $\bLL S \in \aut(\cU/\bLL)$. In the future, we shall abuse slightly of notation and write simply $\Phi\bLL S = \Phi S \bLL$. We say that a toroid $\UoverL$ is a \emph{$k$-orbit toroid} if $\autUoverL$ has $k$ orbits on flags. Following \cite{ARP}, regular toroids are precisely $1$-orbit toroids.

Observe that every translation of $\cU$ induces an automorphism of $\UoverL$ and the translations of $\bLL$ induce the trivial automorphism. Also, the \emph{central inversion of $\E$} $\minid: x \mapsto -x$ is always an automorphism of $\cU$ (see \cite[Table 1]{hubardOrbanicPellicerWeissEquiv4Toroids}) and preserves every lattice $\LL$, so it projects to an automorphism of every toroid. This implies that $ \langle \T(\cU), \minid \rangle \leq \normGU(\bLL)$. Furthermore, since $\minid$ normalizes $\T(\cU)$ and $\T(\cU) \cap \langle \chi \rangle = \{\id\}$ it follows that $\langle \T(\cU), \minid \rangle = \T(\cU) \rtimes \langle \minid \rangle$. Therefore, groups of automorphism of toroids are induced by groups $K$ such that $\T(\cU) \rtimes \langle \minid \rangle \leq K \leq \G(\cU)$. 

By the Correspondence Theorem for groups, those groups $K$ with $\T(\cU) \rtimes \langle \minid \rangle \leq K \leq \G(\cU)$ are in one-to-one correspondence with groups $K'$ such that $ \langle \minid \rangle \leq K' \leq \Go(\cU)$. In this correspondence, the group $K'$ corresponds with the group $\T(\cU) \rtimes K'$.	

Recall that if $\bLL = S^{-1} \bLL' S$ for some $S \in \G(\cU)$ then the toroids $\UoverL$ and $\UoverL'$ are isomorphic. Furthermore, the corresponding automorphism groups are $K/\bLL$ and $(S^{-1} K S) /\bLL'$, for some group $\T(\cU) \rtimes \langle \minid \rangle \leq K \leq \G(\cU)$. Hence, in order to classify toroids up to isomorphism it is sufficient to determine their automorphism groups up to conjugacy. According to the discussion above, we only need to  the find conjugacy classes of groups $K'$ such that $\langle \minid \rangle \leq K' \leq \Go(\cU)$. Furthermore, according to \cite{orbanicPellicerWeiss_MapsOperationsKorbitMaps} the number of flag-orbits of a toroid $\UoverL$ under $\aut(\UoverL)$ is the same as the index of $K$ in $\G(\cU)$ which is the same as the index of $K'$ in $\Go(\cU)$.

We summarize the discussion above in the following lemma. This is essentially Lemma 6 of \cite{hubardOrbanicPellicerWeissEquiv4Toroids}.

\begin{lemma}\label{lem:toroids}
	With the notation above, the following statements hold.
	\begin{enumerate}
		\item A symmetry $S \in \Go(\cU)$ projects to an automorphism of $\UoverL$ if and only if $S$ normalizes $\bLL$. This occurs if and only if $\LL S = \LL$.
		\item Since all lattices are centrally symmetric, $\minid:x \mapsto -x$ always projects to an automorphism of $\UoverL$.
		\item The automorphism group $\autUoverL$ of $\UoverL$ is isomorphic to \[(K' \ltimes \T(\cU)	) / \bLL \cong K' \ltimes (\T(\cU) / \bLL )\] where $K'=\{S \in \Go(\cU) \colon\ S^{-1} \bLL S = \bLL \} =\{S \in \Go(\cU) \colon\  \LL S = \LL \}$. In particular $\langle \minid \rangle \leq K' \leq \Go(\cU)$. The group $\autUoverL$ has $k$ orbits on the set of flags of $\UoverL$ if and only if the index of $K'$ in $\Go(\cU)$ is $k$.
		\item The toroids $\UoverL$ and $\UoverL'$ are isomorphic if and only if $\bLL$ and $\bLL'$ are conjugate in $\G(\cU)$. This in turn is true if and only if exists $S \in \Go(\cU)$ such that $\LL S = \LL'$.
	\end{enumerate}

\end{lemma}

An important class of lattices are those that are preserved by a reflection on a hyperplane. Assume that $\Pi$ is a hyperplane that contains $o$, a lattice $\LL$ is a \emph{vertical translation lattice} (with respect to $\Pi$) if \[\LL= \bigcup_{k \in \bZ} \LL_{0} + ku \] where $\LL_{0} = \LL \cap \Pi$ and $u \in \Pi^{\perp}$. 

Vertical translation lattices with respect to a hyperplane $\Pi$ are trivially preserved by the refection on $\Pi$. However, even those lattices that are not vertical translation are provided with an interesting structure, which is described in the following lemma.

\begin{lemma}\label{lem:refLattices}
	Let $\LL$ be a rank-$n$ lattice on $\E$. Let $R$ be the reflection on a hyperplane $\Pi$ that contains $o$ and assume that $\LL$ is preserved by $R$, that is $\LL R = \LL$. Then the following statements hold:
	\begin{enumerate}
		\item \label{inc:layers} If $\LL_{0} = \LL \cap \Pi$, then $\LL_{0}$ is a rank $n-1$ lattice and \[\LL= \bigcup_{k \in \bZ} \LL_{0} + kw \] for any $w \in \LL \sm \Pi$ with minimum (positive) distance $d$ to $\Pi$.
		\item \label{inc:vert} The point $w$ may be chosen in $\Pi^{\perp}$ if and only if $\LL$ is a vertical translation lattice with respect to $\Pi$.
		\item \label{inc:nonvert} If $\LL$ is not a vertical translation lattice and $\{v_{1}, \dots, v_{n-1}\}$ is a basis for $\LL_{0}$ then $w$ may be chosen of the form  $\frac{1}{2}\left( \alpha_{1} v_{1} + \cdots + \alpha_{n-1} v_{n-1} + u \right)$ where $u \in \Pi^{\perp} \sm \{o\}$, $|u| = 2d$ and $\alpha_{1}, \dots , \alpha_{n-1} \in \{0,1\}$ not all zero. Furthermore, the choice of such $\alpha_{1}, \dots , \alpha_{n-1}$ is unique.
	\end{enumerate}

	\begin{proof}
		To see that $\LL_{0} = \LL \cap \Pi$ is a rank $n-1$ lattice take $x_{1} \in \LL \sm \Pi^{\perp}$ and for $2 \leq i \leq n-1$ take $x_{i} \in \LL$ recursively such that $x_{i}$ does not belong to the $i$-dimensional subspace that contains $x_{1}, \dots , x_{i-1}$ and $\Pi^{\perp}$. Observe that $\{x_{i}+x_{i}R : i \in \{1, \dots , n-1\}\}$ is a set of $n-1$ linearly independent points in $\LL \cap \Pi$. Among all those sets take $\{y_{1}, \dots, y_{n-1}\}$ such that the volume of the parallelotope spanned by $y_{1}, \dots, y_{n-1}$ is minimum, then $\{y_{1}, \dots, y_{n-1}\}$ is a basis for $\LL_{0}$.
		
		Let $w$ be a vector in $\LL \sm \Pi$ with minimum (positive) distance $d$ to $\Pi$. It is clear that $ \bigcup_{k \in \bZ} \LL_{0} + kw \subset \LL$. Assume the other inclusion does not hold. In this situation there must be a point $v$ of $\LL$ in between the hyperplanes $\Pi + kw$ and $\Pi + (k+1)w$ for some $k \in \bZ$. This implies that the point $v-kw$ is a point of $\LL$ with distance strictly less than $d$, which contradicts the choice of $w$. This proves part \ref{inc:layers}.
		
% 		otherwise there would be a point of $\LL$ between the planes $\Pi$ and $\Pi+w_{0}$, contradicting the choice of $w_{0}$. Note also that $u=w_{0}-w_{0}R \in \LL \cap \Pi^{\perp}$. Observe that $u$ has minimum distance $2d$ to $o$ among the points of $\Pi^{\perp} \cap \LL$ since the only possibility would be that there exists a point $u' \in \Pi^{\perp} \cap \LL$ such that $d(u,o) = d$, but this would imply that $\LL$ is a vertical translation lattice. This proves the first part.
		
		If $w$ may be chosen in $\Pi^{\perp}$ then $\LL$ is a vertical translation lattice. Conversely, if $\LL$ is a vertical translation lattice, then the minimality of $d$ implies that $w$ may be chosen in $\Pi^{\perp}$. This proves part \ref{inc:vert}.

		Assume that $\LL$ is not a vertical translation lattice. Let $\{v_{1}, \dots v_{n-1}\}$ be a basis for $\LL_{0}$ and $w \in \LL \sm \Pi$ such that $d=d(w,\Pi)$ is minimum among the points of $\LL \sm \Pi$. Let $u=w-wR$ and note that $|u|=2d$, hence $u$ is a closest point of $\LL \cap \Pi^{\perp}$ to $o$ other than $o$ itself. Observe that $2w-u \in  \LL_{0} $, thus there exist $m_{1}, \dots , m_{n-1} \in \bZ$ such that $2w-u = m_{1}v_{1} +  \cdots + m_{n-1}v_{n-1}$. Let $k_{i} \in \bZ$ for $1 \leq i \leq n-1$ such that $k_{i}$ is even and $\alpha_{i}= m_{i} - k_{i} \in \{0,1\}$. Define \[w_{1} =\frac{1}{2}\left( \alpha_{1}v_{1} +  \cdots + \alpha_{n-1}v_{n-1} + u\right)\] and observe that $\LL= \bigcup_{k \in \bZ} \LL_{0} + kw_{1}$ since $w$ and $w_{1}$ differ by an element in $\LL_{0}$ and hence $d(\Pi,w_{1})=d$. Now, observe that if $\alpha_{i} =0$ for all $1 \leq i \leq n-1$, then all $m_{1}, \dots , m_{n-1}$ are all even and this will imply that $\frac{1}{2}u \in \LL$, contradicting that $\LL$ is not a vertical translation lattice. Finally, suppose there exist $\alpha'_{1}, \dots, \alpha'_{n-1} \in \{0,1\}$ such that if $w_{2} =\frac{1}{2}\left( \alpha'_{1}v_{1} +  \cdots + \alpha'_{n-1}v_{n-1} + u\right)$ satisfies that $\LL = \bigcup_{k \in \bZ} \LL_{0} + kw_{2}$. Let $i \in \{1, \dots, n-1\}$ such that $\alpha_{i} \neq \alpha'_{i}$. On one hand, $w_{1} + w_{2} - u$ is a point in $\LL_{0}$, on the other hand, in the (unique) expression of $w_{1} + w_{2} - u$ as a linear combination of $\{v_{1}, \dots , v_{n-1}\}$, the coefficient of $v_{i}$ is $\frac{1}{2}$, contradicting that $\{v_{1}, \dots , v_{n-1}\}$ is a basis for $\LL_{0}$. \qed
	\end{proof}

\end{lemma}

%discussion about correspondence theorem, conjugacy classes and number of obits.

%lattices%lattices
%equivelar toroids
%symmetries of toroids
%known results about equivelar toroids ¿No es mejor mencionarlos en la introducción?
%result about lattices preserved by an hyperplane reflection
%equivelar toroids
%symmetries of toroids
%known results about equivelar toroids ¿No es mejor mencionarlos en la introducción?
%result about lattices preserved by an hyperplane reflection

%Recall that for the group $\mathbb{E}^n$ a lattice is a discrete subgroup of such group. %ref%
%We define the integer $n$-lattice in $\mathbb{R}^n$ as the points in $\mathbb{E}^n$ with integer coordinates and we define the $n$-cube tiling as the subjacent abstract polytope that arises from said lattice.

%% file: cubic.tex
	%Un poquito de estructura del grupo
	% normales del abeliano 
	% subgrupos de indice 2, todos continen al (C_{2}^{n})^{+} \rtimes A_n
	% Latices invariantes bajo (C_{2}^{n})^{+} \rtimes A_n
	% clasificación de toros de 2 órbitas
	% el único grupo de indice k con 2 < k < n es (C_{2}^{n})^{+} \rtimes A_n
	% clasificación de toros de k órbitas con 2 < k < n 
	% el unico grupo de indice n en S_{n} es S_{n-1} <- Tero cree que ya valió pito y debe estudiar (Dixon)
	% Clasificación de toros de n órbitas. <- Salvo la valida de pito de arriba, ya están. 

In this section we discuss \emph{few-orbit cubic toroids}, meaning equivelar few-orbit $(n+1)$-toroids type $\{4,3^{n-2},4\}$. 
These toroids are quotients of the Euclidean tessellation $\cU = \{4,3^{n-2}, 4\}$ of $\mathbb{E}^{n}$ by unitary $n$-cubes. 

Since the vertex set of $\cU$ is a lattice, according to the discussion in Section\nobreakspace \ref {sec:toroids}, the symmetry group of $\cU$, is the semidirect product $\T(\cU) \rtimes \Go$, where $\T(\cU)$ is the translation group of $\cU$ and $\Go$ is the stabilizer of the base vertex $o$, that is, the group of the vertex figure at $o$ of $\cU$.  

The group $\Go$ is the string Coxeter group $[3^{n-2}, 4]$, generated by reflections $R_{1}, \dots , R_{n}$. The group $\symn = \langle R_{1}, \dots , R_{n-1} \rangle$ is isomorphic to the symmetric group in $n$ symbols, acting on the points of $\E$ by permutation of coordinates. Observe that $R_{n}$ together with its conjugates under $\symn$ generate a group isomorphic to the group $C_{2}^{n}$; we denote this group by $\abeln$. Observe that the vector $(1^{i-1}, -1, 1^{n-i}) \in C_{2}^{n}$ may be identified with the reflection $E_{i}$ in $\abeln$ on the coordinate plane $x_{i} = 0$, given by $(x_{1}, \dots , x_{n}) \mapsto (x_{1}, \dots ,x_{i-1}, -x_{i}, x_{i+1}, \dots, x_{n})$. It is not hard to see that $\Go = \abeln \rtimes \symn$. Finally, observe that the action of conjugacy of $\symn$ on $\abeln$ is precisely the same as the action of $S_{n}$ permuting coordinates in $C_{2}^{n}$. We will use this structure to prove some results regarding the conjugacy classes of subgroups of $\Go$. During the discussion we will identify the elements of $\symn$ with permutations and the elements of $\abeln$ with vectors with entries $\pm 1$ as described before.

\begin{lemma} \label{lem:preservedByAn}
 Let $n \geq 3$ and let $H$ be a subgroup of $\abeln$. Let $\altn$ denote \emph{rotational subgroup} of $\symn$, that is, the group consisting of orientation preserving isometries of $\symn$. If $H$ is normalized by $\altn$ then one of the following hold:
 \begin{enumerate}
  \item $H= \{ \triv \}$.
  \item $H=\langle \minid \rangle$.
  \item $H= \abelnplus$, the rotational subgroup of $\abeln$.
  \item $H= \abeln$.
 \end{enumerate}
 Where $\minid$ denotes the central inversion $x \mapsto -x$.

 \begin{proof}
	We use the structure discussed above for the group $\stabcube$. In this context, the group $\altn$ corresponds to the alternating group and hence we only have to classify groups $H \leq \abeln$ preserved by the action of $\altn$ in the coordinates of its elements.
	
	Suppose that $H \neq \{ \triv \}$ and define $m$ to be the minimum positive number of $-1$ entries in a non-trivial vector of $H$. Let  $A \in H$ be the transformation given by the vector $(a_1,...,a_n)$ and assume that $(a_1,...,a_n)$ has precisely $m$ entries equal to $-1$. 

	If $m=1$, then $A=E_{i}$, the reflection given by the vector that has $1$ in all its coordinates except in the $i$-th. Since $\altn$ acts transitively on the coordinates, for every $j \in \{1, \dots, n\}$ there exists $S_{j} \in \altn$ such that $S_{j}^{-1} E_{i} S_{j} = E_{j}$. Since $H$ is preserved by $\altn$ and the set $\{ E_{j}: 1 \leq j \leq n \}$ generates $\abeln$, $H$ must be $\abeln$.
	
	If $m = 2$ then $A=E_{i} E_{j}$, the transformation given by vector whose $i$-th and $j$-th coordinates are $-1$. Proceeding in a similar way as before we may show that $\{E_{i}E_{j} : 1 \leq i < j \leq n\} \subset H$, which implies that $\abelnplus \leq H$. However, by the minimality of $m$ we must have that $H = \abelnplus$.
	
	If $2<m<n$, without lose of generality we may assume that $a_i= -1$ if and only if $1\leq i\leq m$, that is  $A$ is the transformation induced by the vector $(-1^{m}, 1^{n-m})$. Since $H$ is preserved by $\altn$, the transformation $A'$ given by $(1,-1^{m}, 1^{n-m-1})$ belongs to $H$. Therefore, $A A' \in H$ is given by the vector $(-1, 1^{m-1}, -1, 1^{n-m-1} )$, which contradicts the minimality of $m$.
	
	Finally, if $m=n$ then $H = \langle \minid \rangle$. \qed
\end{proof}
\end{lemma}

Recall that if a group $G$ is the semidirect product $N \rtimes H$ there exists a natural mapping $\eta :G \to H$ whose kernel is $N$. Furthermore, if $K \leq G$ then the restriction of $\eta$ to $K$ has kernel $K \cap N$; in particular if $G$ is finite we have 
\begin{equation}\label{eq:indexIm}
	[G:K] = \frac{|G|}{|K|} = \frac{|N| |H|}{|K|} \geq  \frac{|N\cap K| |H|}{|K|} = \frac{|H|}{|\eta(K)|} = [H:\eta(K)]
\end{equation}

\subsection{Cubic toroids with \texorpdfstring{$2$}{2} orbits.} 

In this section we classify $2$-orbit toroids of type $\{4, 3^{n-2}, 4 \}$. According to Lemma\nobreakspace \ref {lem:toroids}, to do so we need to find lattices preserved by index two subgroups of $[4, 3^{n-2}, 4]$. The following results enumerate such groups.

\begin{lemma}\label{lem:index2Cube}
 Let $n \geq 3$. If $K\leq \stabcube$ is an index two subgroup, then $\abelnplus \rtimes \altn \leq K$. 
 \begin{proof}
  Let $\eta: \stabcube \to \symn$ be the natural mapping. By Equation\nobreakspace \textup {(\ref {eq:indexIm})}, $[\symn:\eta(K)] \leq 2$, therefore $\altn \leq \eta(K)$. This implies that $|K\cap \abeln| \geq 2^{n-1}$. Since $K$ and $\abeln$ are normal subgroups, $K \cap \abeln$ must be preserved by conjugation under $\altn$, and by Lemma\nobreakspace \ref {lem:preservedByAn}, $\abelnplus \leq K$.
  
  Since $\altn \leq \eta(K)$, for every $3$-cycle $S$ in $\altn$ there exists $A \in \abeln$ such that $AS \in K$. This implies that \[(AS)^{2} = ASAS = ASAS^{-1}  S^{2}  \in K.\] However, $ASAS^{-1} \in \abelnplus \leq K$ which implies that $S^{2} = S^{-1}\in K$. Since this holds for every $3$-cycle $S$, then $\altn \leq K$ which implies that $\abelnplus \rtimes \altn \leq K$, as desired. \qed
 \end{proof}
\end{lemma}

\begin{corollary}\label{cor:index2Cube}
 If $K$ is an index two subgroup of $\stabcube$, then $K$ is one of the following.
 \begin{enumerate}
  \item $(\stabcube)^{+}$, the rotational subgroup of $\stabcube$.
  \item $\abelnplus \rtimes \symn$.
  \item $\abeln \rtimes \altn$.
 \end{enumerate}
	\begin{proof}
	 It is clear that those three groups are different. By Lemma\nobreakspace \ref {lem:index2Cube}, $K$ must contain $\abelnplus \rtimes \altn$ which is a normal subgroup of $\stabcube$ of index $4$. By the Correspondence Theorem for groups there must be at most three index-two subgroups containing $\abelnplus \rtimes \altn$. \qed
	\end{proof} 
\end{corollary}

Corollary\nobreakspace \ref {cor:index2Cube} determines all the subgroups of $\stabcube$ of index $2$. According to Lemma\nobreakspace \ref {lem:toroids}, by classifying the lattices preserved by such groups we obtain a classification of $2$-orbits toroids.

By Lemma\nobreakspace \ref {lem:index2Cube}, every lattice preserved by an index 2 subgroup of $\stabcube$ is also invariant under $\abelnplus\rtimes\altn$. Therefore it is useful to know all those lattices.

Consider the vectors ${v_1}:=(1,1, \dots, 1)$, ${v_k}:=( -1,1^{k-2} ,-1 , 1^{n-k})$ for $2 \leq k \leq n$ and ${w_k}:=(1^{k-1}, -1 , 1^{n-k})$ for $1\leq k \leq n$. Let $\latEven$ and $\latOdd$ be the lattices whose basis are $\{v_{i} : 1 \leq i \leq n \}$ and $\{w_{i} : 1 \leq i \leq n \}$, respectively. First observe that for every $s \in \bN$, we have that \begin{equation}\label{eq:latEvenlatOdd}
4s \cln \subset 2s \fcln \subset s(\LL_{0} \cap \LL_{1}).
\end{equation}
In particular, $ 2e_{j} + 2e_{k} \in \LL_{0}\cap \LL_{1}$, for $j,k \in \{1, \dots, n\}$, $j \neq k$. This implies that $\LL_{0}$ and $\LL_{1}$ are preserved by $\abelnplus$, since $ v_{k}-v_{k}(E_{i}E_{j}) \in 2\fcln$ and $w_{k}-w_{k}(E_{i}E_{j}) \in 2\fcln$ for $i,j,k \in \{1, \dots, n\}$, $i < j$. Similar arguments prove that $\LL_{0}$ and $\LL_{1}$ are also preserved by $\altn$. Now it is easy to classify the lattices preserved by $\abelnplus \rtimes \altn$.

% \begin{lem}\label{lem:latEvenOdd}
%  Let $\LL$ be a rank-$n$ lattice preserved by $\abelnplus \rtimes \symn$
% \begin{proof}
%   Let $\bar{e}_i:=(0^{i-1},1,0^{n-i})$ for $1\leq i \leq n$, and notice that $\mathbf{a_j}-\mathbf{a_k}=2\bar{e}_j+2\bar{e}_k\in \latEven$ for $1 \leq i<j\leq n$. Therefore, $\mathbf{a_1} \sigma, \ldots, \mathbf{a_n}\sigma \in \latEven$ for every $\sigma \in \abelnplus$. Furthermore, since, for the transposition $(1\,j)\in\symn$, $\mathbf{a_i}(1\,j)=\mathbf{a_1}-2\bar{e}_i-2\bar{e}_j\in\latEven$, it follows that $\latEven$ is preserved by the action of $\symn$. We conclude that $\latEven$ is invariant under $\abelnplus \rtimes \symn$.
%   
%   For $1\leq i<j \leq n$, we have $\mathbf{b_i} \pm \mathbf{b_j} = 2\bar{e}_i\pm2\bar{e}_j\in \latOdd$, therefore $\mathbf {b_i}\sigma \in \latOdd$ for every $\sigma\in\abelnplus$ and $1\leq i \leq n$. And since it is immediate that $\symn$ preserves $\{\mathbf{b_1},\ldots,\mathbf{b_n}\}$, it follows that $\latOdd$ is invariant under $\abelnplus \rtimes \symn$.
% \end{proof}
% \end{lem}

\begin{lemma}\label{lem:laticesC2plusAlt}
 If $n\geq 4$ and is $\LL$ a rank $n$ integer lattice preserved by the group $\abelnplus \rtimes \altn$, then $\LL$ is an integer multiple of one of the following:
 \begin{enumerate}
  \item $\regLat$ for $k \in \{1,2,n\}$.
  \item $\latEven$.
  \item $\latOdd$.
 \end{enumerate}
\begin{proof}
 Let $s$ be the smallest positive integer among the coordinates of the vectors of $\LL$. As before, since $\altn$ acts transitively on the coordinates of $\E$, $s$ divide every entry of every point of $\LL$. By permuting coordinates we may assume that $(s,s_{2}, \dots ,s_{n}) \in \LL$ for some $s_{2}, \dots, s_{n} \in \bZ$. By adding some elements in $\LL$ and using elements of $\abelnplus \rtimes \altn$ we deduce the existence of some vectors in $\LL$ as follows
 \begin{align*}
    (s,s_{2}, \dots ,s_{n}) \in \LL & \Rightarrow (s,s_2,\ldots,s_n)E_{1,2}=(-s,-s_{2}, s_{3} \dots ,s_{n}) \in \LL, \\     
    &\Rightarrow (s,s_{2}, \dots ,s_{n})-(-s,-s_{2}, \dots ,s_{n}) = \\
	  &\hphantom{\Rightarrow} = (2s,2s_{2}, 0, \dots ,0) \in \LL, \\
    &\Rightarrow (2s,2s_{2}, 0, \dots ,0)E_{2,3}=(2s,-2s_{2},-0, \dots ,0) \in \LL, \\
    &\Rightarrow (2s,2s_{2}, 0, \dots ,0)+(2s,-2s_{2}, 0, \dots ,0)= \\ 
    &\hphantom{\Rightarrow} =(4s,0, \dots ,0) \in \LL, \\
    &\Rightarrow (0^{j-1},4s, 0^{n-j}) \in \LL \text{ for all } 1 \leq j \leq n.
 \end{align*}
  This implies that $4s\cln \subset \LL $. Therefore we may write every point of $\LL$ as $(s_{1}, \dots, s_{n}) + x$ where $s_{1}, \dots, s_{n} \in \{-s,0,s,2s\}$ and $x \in 4s\cln$.
  
  Assume there exists a point of the form $(s, s_{2}, \dots, s_{n}) + x$ with $s_{2}, \dots, s_{n} \in \{-s,0,s,2s\}$ and $x \in 4s\cln$, such that for some $i \in \{2, \dots, n\}$, $s_{i} \in \{0,2\}$. By permuting coordinates, we may assume that $i=2$. As seen before, $(2s,2s_{2}, 0, \dots ,0) \in \LL$, but $4s|2s_{2}$ and since $(0, 4s, 0, \dots,0) \in \LL$, then  $(2s, 0, \dots ,0) \in \LL$ and hence $(0^{j-1},2s, 0^{n-j}) \in \LL$ for all $ 1 \leq j \leq n$, which implies that $\LL$ is preserved not only by $\abelnplus \rtimes \altn$ but by $\abeln \rtimes \altn$. This implies that $\LL$ is $s \regLat$ with $k \in \{1,2,n\}$ (see \cite[p. 166]{ARP}).
  
  If every point of $\LL$ is of the form $(\pm s, \dots ,\pm s) + x$ with $x \in 4s \cln$ there are several cases to consider. Assume there exist points $p,q \in \LL$ with $p = (p_{1}, \dots, p_{n}) + x_{p}$, $q = (q_{1}, \dots, q_{n}) + x_{q}$ where $|p_{i}|= |q_{i}| = s $ for $i \in \{1, \dots, n\}$, $x_{p}, x_{q} \in 4s\cln$ such that the number of entries of $(p_{1}, \dots, p_{n})$ equal to $-s$ is even and the number of entries of $(q_{1}, \dots, q_{n})$ equal to $-s$ is odd. By performing some permutation of coordinates and even number of sign changes (if needed) we may assume that $p= (-s^{k}, s^{n-k}) +x_{p}$ and $q= (-s^{m}, s^{n-m}) +x_{q}$ with $k$ even, $m$ odd and $k < m$. In this situation $p-q = (0^{k}, 2^{m-k}, 0^{n-m})+(x_{p}-x_{q})$ and since $m-k$ is odd and $2se_{i} + 2se_{j} \in \LL$ for every $i,j \in \{0, \dots, n\}$, then $2se_{i} \in \LL$ for every $i \in \{1, \dots n\}$. This implies that $\LL$ is preserved by $\abeln \rtimes \altn$. As before, this implies that $\LL$ is $s\regLat$ for $k \in \{1,2,n\}$. 
  
  Now, with the notation used above, the  other possibility is that for every point $p = (p_{1}, \dots, p_{n}) + x_{p}$, the parity of the number of entries of $(p_{1}, \dots, p_{n})$ equal to $-s$ does not depend on $p$. Recall that $2s \fcln \subset \latEven \cap \latOdd$. Hence, if the parity is even then all such points $(p_{1}, \dots, p_{n})$ belong to $s\latEven$, if the parity is odd then all belong to $s \latOdd$. Furthermore, since $4s \cln \subset s\latEven \cap s\latOdd$, then $\LL \subset s\LL_{i}$, for exactly one $i\in \{0,1\}$. The other inclusion follows from the fact that $2s(e_{i} +e_{j}) \in \LL $ for every $i,j \in \{1, \dots ,n\}$, which implies that $\LL$ contains either $\{sv_{i} :  1 \leq i \leq n\}$ or $\{sw_{i} : 1 \leq i \leq n\}$. \qed
\end{proof}
\end{lemma}

We say that a $2$-orbit toroid belongs to \emph{class $2_{I}$} for $I \subseteq \{0, \dots ,n\}$ if given a flag $\Phi$ and $i \in I$, then the $i$-adjacent flag $\Phi^{i}$ belongs to the same orbit as $\Phi$. In particular, chiral toroids are those in class $2_{\emptyset}$.

The following result follows almost immediately from Lemma\nobreakspace \ref {lem:laticesC2plusAlt} and together with \cite[Theorem 7]{hubardOrbanicPellicerWeissEquiv4Toroids} completes a classification of $2$-orbit toroids of type $\{4, 3^{n-2}, 4\}$.

\begin{theorem}\label{thm:cubic2orbits}
 Let $n \geq 3$. If $n$ is odd, there are no equivelar $(n+1)$-toroids of type $\{4, 3^{n-2}, 4\}$ with $2$ flag-orbits. If $n$ is even, every equivelar $(n+1)$-toroid of type $\{4, 3^{n-2}, 4\}$ with two flag-orbits may be described as a toroid $\{4, 3^{n-2}, 4\}_{\bLL}$ with $\bLL$ an integer multiple of the lattice group $\mathbf{\latOdd}$ described above. Furthermore, such toroids belong to class $2_{\{1, 2, \dots, n-1\}}$.
 \begin{proof} 
  According to Lemma\nobreakspace \ref {lem:toroids}, every $2$-orbit cubic toroid must be given by a lattice $\LL$ preserved by an index two subgroup $K$ of $\Go$ containing $\minid$. By Lemma\nobreakspace \ref {lem:index2Cube}, the only possibilities for $K$ are $(\stabcube)^{+}$, $\abeln \rtimes \altn$ and $\abelnplus \rtimes \symn$. However, $K$ cannot be $(\stabcube)^{+}$, since this will induce a chiral $(n+1)$-toroid and they are known not to exist if $n \geq 3$ (see \cite[Theorem 9.1]{mcMullenShulte_HigherToroidalRP} and \cite[Section 6D]{ARP}). As mentioned before, if $\LL$ is a lattice preserved by $\abeln \rtimes \altn$ then $\LL$ is an integer multiple of a lattice $\regLat$ for $k \in \{1,2,n\}$, but such lattices induce regular toroids. Therefore, the only possibility for $K$ is $\abelnplus \rtimes \symn$. 
  
  Note that if $n$ is odd, then $\minid \notin \abelnplus \rtimes \symn$. This implies that there are no two-orbit $(n+1)$-toroids whenever $n$ is odd. If $n$ is even then $\LL$ must be preserved by $\abelnplus \rtimes \symn$. In particular, $\LL$ must be preserved by $\abelnplus \rtimes \altn$ and those lattices are classified in Lemma\nobreakspace \ref {lem:laticesC2plusAlt}. Since the lattices $s \regLat$ induce regular toroids for all $s \in \bZ$, then $\LL$ must be an integer multiple of $\latEven$ or an integer multiple of $\latOdd$. Finally, observe that both $\latEven$ and $\latOdd$ are preserved by $\symn$ and since they are isometric lattices, then they produce isomorphic toroids. Consequently, every two-orbit toroid is induced by an integer multiple of $\latOdd$, as stated.
  
  To see that those toroids belong to class $2_{\{1, 2, \dots, n-1\}}$ note that if $R_{0}, \dots , R_{n}$ are the distinguished generators of $\G(\cU)$ then $R_{1}, \dots, R_{n-1} \in \symn$ and hence preserve $\latOdd$. These generators induce automorphisms $\bar{R}_{1}, \dots, \bar{R}_{n-1}$ of $\cU / \bLL$ such that if $\Phi \bLL$ is the base flag of $\cU / \bLL$, then $\Phi \bLL \bar{R}_{i} = \Phi R_{i} \bLL = \Phi^{i} \bLL = (\Phi \bLL)^{i}$, for $1 \leq i \leq n-1$. Moreover, there is no automorphism of $\cU/ \bLL$ mapping $\Phi \bLL$ to $(\Phi \bLL)^{i}$ for $i \in \{0,n\}$ since this will imply that $\LL$ is preserved by $\abeln \rtimes \altn$. Therefore, $\UoverL$ belongs to class $2_{\{1,2, \dots, n\}}$. \qed
  
 \end{proof} 
\end{theorem}

\subsection{Cubic \texorpdfstring{$(n+1)$}{(n+1)}-toroids with \texorpdfstring{$k$}{k} flag-orbits for \texorpdfstring{$2<k<n$}{2<k<n}.}
Now we proceed to classify toroids of type $\{4,3^{n-2},4\}$ with $k$ orbits for $2 < k < n$. The following lemma is the key of such classification.

\begin{lemma}\label{lem:cubicNoIndexkGroups}
 Let $n \geq 5$. If $K \leq \stabcube$ has index $k$ for some $2 < k < n$ then $K = \abelnplus \rtimes \altn$. 
 \begin{proof}
  Assume $K$ as above and let $\eta: \stabcube \to \symn$ be the projection to $\symn$. By Equation\nobreakspace \textup {(\ref {eq:indexIm})}, $[\symn : \eta(K)] \leq k$ and since $n \geq 5$, $\eta(K)$ must contain $\altn$. Let $K_{0}=K \cap \abeln$. Observe that $K_{0}$ is normalized by $\altn$ and those groups are classified in Lemma\nobreakspace \ref {lem:preservedByAn}. If $\abelnplus \leq K_{0}$, then we may proceed as in the proof of Lemma\nobreakspace \ref {lem:index2Cube} and conclude that $\abelnplus \rtimes \altn \subset K$ and $k \leq 4$, which implies that $\abelnplus \rtimes \altn = K$. 
  If $K_{0} \in \{\{ \triv \}, \langle \minid \rangle \}$ then we have that $|K| = |K_{0}||\eta(K)| \leq 2 n!$, hence \[n > [\stabcube:K] = \frac{2^{n}n!}{|K|} \geq \frac{2^{n}n!}{2n!} = 2^{n-1}\] which is a contradiction since $n \geq 5$. \qed
 \end{proof}
\end{lemma}

As an immediate corollary we have the following result.

\begin{theorem}\label{thm:nokOrbitsToroids}
 Let $n \geq 5$. Then there are no equivelar $(n+1)$-toroids of type $\{4, 3^{n-2},4\}$ with $k$ flag-orbits whenever $2 < k < n$.
%  \begin{proof}
%   According to \cref{lem:cubicNoIndexkGroups}, a cubic $(n+1)$-toroid with $k$ orbits, for $2 < k < n$ must be induced by a lattice preserved by $\abelnplus \rtimes \altn$. However such lattices are classified in Lemma\nobreakspace \ref {lem:laticesC2plusAlt} and they induce either regular toroids or two-orbit toroids. 
%  \end{proof}
\end{theorem}

So far we have classified cubic two-orbit $(n+1)$-toroids for arbitrary $n \geq 3$. We also proved the non-existence of cubic $k$-orbit $(n+1)$-toroids for $2 < k < n$ when $n \geq 5$. These results together with those in \cite[Section 6D]{ARP} and  \cite{hubardOrbanicPellicerWeissEquiv4Toroids} almost complete the classification of equivelar few-orbit toroids of type $\{4, 3^{n-1}, 4\}$. In order to complete this classification we need to classify four dimensional toroids with three flag-orbits and $n$-dimensional toroids with $n$ flag-orbits. 

The following results classify cubic $(4+1)$-toroids with three flag-orbits. In order to give such classification first consider the group $\died \leq \symn[4]$ generated by the reflections $(x_{1}, x_{2},x_{3},x_{4}) \mapsto (x_{2}, x_{1},x_{4},x_{3}) $ and $(x_{1}, x_{2},x_{3},x_{4}) \mapsto (x_{1}, x_{4},x_{3},x_{2})$. This group acts in the coordinates of the points of $\E[4]$ as the dihedral group acts on the vertices of an square labeled by $1$, $2$, $3$, $4$.

\begin{lemma}\label{lem:4cubicIndex3Groups}
 Let $\stabcube[4]$ be the stabilizer of the origin in the cubic tessellation of $\E[4]$. Up to conjugacy, the only subgroup of index $3$ of $\stabcube[4]$ is $\abeln[4] \rtimes \died$.
 \begin{proof}
  Let $K$ be a subgroup of $\stabcube$ of index $3$. Since $K$ is a 2-Sylow subgroup of $\stabcube$ and $\abeln[4]$ is a normal 2-subgroup, $\abeln[4]$ must be contained in $K$ and hence $K = \abeln[4] \rtimes \eta(K)$, where $\eta$ is the natural mapping to $\symn[4]$. Furthermore, $\eta(K)$ must have index $3$ in $\symn[4]$ which implies that, up to conjugacy, $\eta(K)$ is $\died$. \qed
 \end{proof}
\end{lemma}

\begin{lemma}\label{lem:4cubic3orbitlattices}
 Assume $\LL$ is an integer lattice preserved by the group $\abeln[4] \rtimes \died$, then $\bLL$ is a integer multiple of one of the following:
 \begin{itemize}
  \item $\LL_{(1^{k},0^{4-k})}$ for $k \in \{1,2,4\}$.
  \item $\LatProd$, where $\LatProd$ is the lattice generated by the vectors $(1,0,1,0)$, $(1,0,-1,0)$, $(0,1,0,1)$ and $(0,1,0,-1)$.
 \end{itemize}
	\begin{proof}
	 It is clear that all those lattices are preserved by $\abeln[4] \rtimes \died$. Assume that $\LL$ is a lattice preserved by $\abeln[4] \rtimes \died$ and let $s$ be the minimum positive value among all the coordinates of the points in $\LL$ and take $(s, s_{2}, s_{3}, s_{4}) \in \LL$. Since $\LL$ is preserved by $\abeln[4]$ and $\died$ acts transitively on the coordinates of the elements in $\LL$, proceeding as before we may conclude that $(0^{j-1},2s, 0^{n-j}) \in \LL$ for all $ 1 \leq j \leq n$. Therefore, we may assume that $s_{i} \in \{0,s\}$ for $i \in \{2,3,4\}$. 
	 
	 Among all such points, we take $x \in \LL$ with minimum number $k$ of non-zero entries. If $k = 1$ then $\LL = s \clFour$. If $k=2$ up to permutation of coordinates with an element of $\died$, there are two possibilities: $x = (s,s,0,0)$ in which case $\LL = s \fclFour$; or $ x = (s,0,s,0)$ which implies that $\LL = \LatProd$. The case $k=3$ is impossible since up to the action of $\died$, the only possibility is $x = (s,s,s,0)$ which implies $(0,s,s,s) \in \LL$ and hence $(s,0,0,-s) \in \LL$, contradicting the choice of $k$. The only remaining possibility is $k = 4$ which implies that $\LL = s \bclFour$. \qed
	\end{proof}

\end{lemma}

\begin{theorem}\label{thm:5cubicToroids3orbits}
 Every equivelar rank-$5$ toroid of type $\{4,3,3,4\}$ with three flag-orbits may be described as a toroid $\{4,3,3,4\}_{s \bLL}$, where $s \in \bZ$ and $\bLL$ is the lattice group generated by the translations with respect the vectors $(1,0,1,0)$, $(1,0,-1,0)$, $(0,1,0,1)$ and $(0,1,0,-1)$ . 
	\begin{proof}
	 It follows from Lemma\nobreakspace \ref {lem:4cubic3orbitlattices}, since the integer multiples of the lattices $\LL_{(1,0^{4-k})}$ for $k \in \{1,2,4\}$ induce regular toroids. \qed
	\end{proof}
\end{theorem}

\subsection{Cubic \texorpdfstring{$(n+1)$}{(n+1)}-toroids with \texorpdfstring{$n$}{n} flag-orbits.}

Now we proceed to classify $(n+1)$-cubic toroids with $n$ flag-orbits for $n \geq 4$. The strategy is essentially the same that we have used before: first we determine the $n$-index subgroups of $\stabcube$ and then we determine the lattices preserved by each of them. In order to classify subgroups of $\stabcube$ of index $n$ we use the following result regarding permutation groups. This is a consequence of \cite[Theorem 5.2B]{dixonPermutationGroups}.

\begin{theorem}\label{thm:indexnSubgpsSn}
 Let $n \geq 5$ and let $S_{n}$ be the group of permutations of $\{1, \dots , n \}$. If $G$ is a subgroup of $S_{n}$ of index $n$, then up to conjugacy, one of the following holds:
 \begin{itemize}
  \item The group $G$ is the stabilizer of $\{n\}$ and hence isomorphic to $S_{n-1}$.
  \item $n = 6$ and $G$ acts on $\{1, 2, \dots 6\}$ as the group $PGL_{2}(5)$ acts on the points of the projective line over $GF(5)$, the field with $5$ elements. 
 \end{itemize}
\end{theorem}

Now we may classify the subgroups of $\stabcube$ of index $n$.

\begin{lemma}\label{lem:cubicIndexnSubgps}
 Let $n \geq 4$. If $K$ is a subgroup of $\stabcube$ of index $n$, then up to conjugacy, one of the following holds:
 \begin{itemize}
  \item $K = \abeln \rtimes \symn[n-1]$.
  \item $n = 4$ and $K = \abelnplus[4] \rtimes \altn[4]$.
  \item $n=6$ and $K=\abeln[6] \rtimes \pgl$.
 \end{itemize}
	Here $\symn[n-1]$ denotes the stabilizer of the last coordinate in $\symn$ and $\pgl$ denotes the subgroup of $\symn[6]$ that acts on the coordinates of $\E[6]$ as $PGL_{2}(5)$.
 \begin{proof}
  The result is a consequence of Theorem\nobreakspace \ref {thm:indexnSubgpsSn}. By Equation\nobreakspace \textup {(\ref {eq:indexIm})}, $[\symn:\eta(K)] \leq n$, which implies that $\eta(K)$ is $\symn$, $\altn$, or one of the groups described in Theorem\nobreakspace \ref {thm:indexnSubgpsSn}. Let $K_{0} = K \cap \abeln$ and recall that $|K|= |\eta(K)||K_{0}|$. 
  
  If $\eta(K) = \symn$, then $|K| = (n-1)! 2^{n} = n! |K_{0}|$, which implies that $\frac{2^{n}}{n}= |K_{0}|$. However, if $n \geq 4$ we have $2 < \frac{2^{n}}{n} < 2^{n-1}$ but this is impossible since $K_{0}$ is a subgroup of $\abeln$ normalized by $\symn$ and by Lemma\nobreakspace \ref {lem:preservedByAn}, must be one among $\abeln$, $\abelnplus$ and $\minid$.
  
  Proceeding in a similar way we may see that if $\eta(K) = \altn$, the unique possibility for $K_{0}$ is $\abelnplus$ and this is only possible if $n=4$. With similar arguments to those used in the proof of Lemma\nobreakspace \ref {lem:index2Cube} we may conclude that $K = \altn[4] \rtimes \abelnplus[4]$.
  
  Finally, if $\eta(K)$ has index $n$ in $\symn$, then we have $|K| = (n-1)! 2^{n} = |\eta(K)| |K_{0}|=(n-1)! |K_{0}|$, therefore $K_{0}$ must be $\abeln$ and $K = \abeln \rtimes \symn[n-1]$ or $n=6$ and $K= \abeln \rtimes \pgl$. \qed
 \end{proof}
\end{lemma}

Now we proceed to classify lattices preserved by subgroups of $\stabcube$ of index $n$. Lattices preserved by $\abelnplus[4] \rtimes \altn[4]$ are described in Lemma\nobreakspace \ref {lem:laticesC2plusAlt}. It only remains to determine those rank-$n$ lattices preserved by $\abeln \rtimes \symn[n-1]$ and those rank-$6$ lattices preserved by $\abeln[6] \rtimes \pgl$. Such lattices are described in the following results.

\begin{lemma}\label{lem:latticesPGL25}
	If $\LL$ is an integer rank-$6$ lattice preserved by $\abeln[6] \rtimes \pgl$, then $\LL = s \regLat[6]$ for some $s \in \bZ$ and $k \in \{1,2,6\}$. 
	\begin{proof}
		Let $\LL$ be a rank-$6$ lattices preserved by $\abeln[6] \rtimes \pgl$. Let $s \in \bZ$ be a positive integer such that every entry of any vector in $\LL$ has absolute value at least $s$. As before, since $\pgl$ acts transitively on the coordinates of $\E[6]$, all the entries of any vector of $\LL$ must be multiples of $s$. Let $(s_{1}, \dots , s_{6}) \in \LL$. Without loss of generality, we may assume that $s_{1} = s$. Since $\LL$ is preserved by $\abeln[6]$ this implies that $(-s,s_{2},\dots,s_{6}) \in \LL$. Hence $2s e_{1} \in \LL$ and therefore $2s e_{i} \in \LL$ for $i \in \{1, \dots , 6\}$. 
		
		We may take $(s_{1}, \dots , s_{6}) \in \LL$ such that all its entries are either $0$ or $s$. Furthermore, since the action of $\pgl$ is $3$-transitive on the coordinates of $\E[6]$, we may assume that $(s_{1}, \dots , s_{6})$ is of the form $(s^{k},0^{6-k})$. Among all non-zero points in $\LL$ with this form take one of those where $k$ is minimum.  
		
		If $k =1$ then $se_{i} \in \LL$ for $i \in \{1, \dots,6\}$ which implies that $s\cln[5] \subset \LL$. However, since all the vectors of $\LL$ have as entries integer multiples of $s$, we must have $\LL = s\cln[5]$.
		
		Suppose $k=2$. Since $\LL$ is preserved by $\pgl$, then $s\fcln[4] \subset \LL$. On the other hand, by the observations made above every point of $\LL$ may be written as the sum of  a point in $s\fcln[4]$ and a point in $2s\cln[5]$. Since  $2s\cln[5] \leq s\fcln[4]$ we must have $\LL = s\fcln[4]$. 
		
		Suppose $k \in \{3,4,5\}$ and take $(s^{k}, 0 ^{6-k}) \in \LL$. Given that $\LL$ is preserved by $\pgl$, we must have that $(0, s^{k}, 0^{6-k-1})$ is a point of $\LL$. This implies that $(s,0^{k-1},s,0^{6-k-1}) \in \LL$, hence $(s,s,0,0,0,0) \in \LL$ which contradicts the choice of $k$.
		
		Finally, if $k = 6$, proceding in a similar way of that used when $k=2$, we may conclude that $\LL = s \bcln[6]$. \qed
	\end{proof}
\end{lemma}

The previous lemma implies that all $(6+1)$-toroids induced by lattices preserved by $\abeln[6] \rtimes \pgl$ are regular toroids. Since all $(4+1)$-toroids induced by lattices preserved by $\abeln[4] \rtimes \altn[4]$ are regular or $2$-orbit, Lemma\nobreakspace \ref {lem:cubicIndexnSubgps} and the fact that $\abeln \rtimes \symn[n-1]$ is maximal in $\abeln \rtimes \symn$ imply that if there exist $(n+1)$-toroids with $n$ flag-orbits, they must be induced by lattices preserved by $\abeln \rtimes \symn[n-1]$ that do not induce regular toroids. We now proceed to classify those lattices.

We use Lemma\nobreakspace \ref {lem:refLattices} to classify lattices preserved by $\abeln \rtimes \symn[n-1]$. We shall assume that $\symn[n-1]$ is precisely the stabilizer of the last coordinate on $\E$. Let $\Pi$ be the hyperplane $x_{n}=0$ and $R$ the reflection through $\Pi$, this is $R: (x_{1}, \dots , x_{n}) \mapsto (x_{1}, \dots,x_{n-1}, -x_{n})$. If $\LL$ is a lattice preserved by $\abeln \rtimes \symn[n-1]$, then $\LL$ is a lattice preserved by $R$. Since $\symn[n-1]$ preserves $\Pi$, if $\LL_{0}=\LL \cap \Pi$, then $\LL_{0}$ must be a rank-$(n-1)$ lattice preserved by the restriction of the action of $\abeln \rtimes \symn[n-1]$ to $\Pi$ and therefore $\LL_{0}$ must be one among $s\cln[n-2]$, $s\fcln[n-3]$ and $s\bcln[n-1]$ for some positive integer $s$. 

If $\LL$ is a vertical translation lattice then $\LL = \trasLat[(de_{n})] $ for some $d \in \bN$. All the lattices of this form are preserved by $\abeln \rtimes \symn[n-1]$, we only need to determine those values for $d$ which induce regular lattices. If $\LL_{0}= s\cln[n-2]$, then $d=s$ implies that $\LL = s\cln$. If $\LL_{0}$ is $s\fcln[n-3]$ or $s\bcln[n-1]$, then any value of $d$ gives a lattice that induce an $n$-orbit toroid. Therefore, for the following discussion we shall asume that $\LL$ is not a vertical translation lattice.

Assume that $\LL_{0} = s\cln[n-2]$. Let $u$ be the point of $\LL \cap \Pi^{\perp}$ closets to $o$, say $u= d e_{n}$.  By Lemma\nobreakspace \ref {lem:refLattices} there exist unique $\alpha_{1}, \dots, \alpha_{n-1} \in \{0,1\}$, not all zero, such that if $w=\frac{s}{2}\left(\alpha_{1}e_{1} + \cdots + \alpha_{n-1}e_{n-1}\right) + \frac{d}{2} e_{n} $, then $\LL = \trasLat$. Since $\LL$ is preserved by $\symn[n-1]$ the only possible choice of the numbers $\alpha_{1}, \dots, \alpha_{n-1}$ is $\alpha_{1}= \cdots= \alpha_{n-1}=1$. Since $w$ must be a point with integer coordinates, $s$ and $d$ must be even. If $d={s}$, then $\LL = \frac{s}{2} \bcln$. In any other case $\LL$ induce an $n$-orbit toroid. 

Now assume that $\LL_{0} = s \fcln[n-3]$. As before, let $u$ be  the closest point of $\Pi^{\perp}$ to $o$, say $u=de_{n}$ for some integer $d$. Recall that $\{2se_{1}, s(e_{2}-e_{1}), s(e_{3}-e_{2}), \dots , s(e_{n-1}-e_{n-2})\}$ is a basis for $\LL_{0}$. Take $\alpha_{1}, \dots , \alpha_{n-1} \in \{0,1\}$, not all zero, such that 
\begin{align*}
	w &= \frac{s}{2} \Big( \alpha_{1} (2e_{1}) + \alpha_{2} (e_{2}-e_{1}) + \dots \alpha_{n-1}(e_{n-1}-e_{n-2}) \Big) + \frac{d}{2}e_{n}\\
	&=  \frac{s}{2} \Big( (2\alpha_{1}-\alpha_{2})e_{1} + (\alpha_{2}-\alpha_{3})e_{2}+ \cdots \\
	&\phantom{=} + (\alpha_{n-2}-\alpha_{n-1})e_{n-2} + \alpha_{n-1} e_{n-1}\Big)+ \frac{d}{2}e_{n}
\end{align*}
satisfies that $\LL = \trasLat $. Let $w_{i}$ be the image of $w$ under the reflection on the plane $x_{i}=0$ for $1 \leq i \leq n-1$. Since $\LL$ is preserved by $\abeln$, $w_{i} \in \LL$ for all $i \in \{1, \dots, n\}$. If $\alpha_{n-1}=1$ then $w - w_{n-1} = se_{n-1}$. This contradicts that $\LL_{0} = s \fcln[n-3]$; then $\alpha_{n-1} = 0$. Similarly if $\alpha_{n-1} = 0$ and $\alpha_{n-2} = 1$, then $w - w_{n-1} = se_{n-2}$, hence $\alpha_{n-2} =0$. Proceding in a similar way, we may conclude that $\alpha_{2} = \cdots = \alpha_{n-1}=0$, which implies that $\alpha_{1}=1$ and $w= s e_{1} + \frac{d}{2} e_{n}$. Again, since $w$ must have integer coordinates, then $d$ must be even. If $d=2s$, then $\LL = s\fcln$. In any other case $\LL$ is a lattice that induce an $n$-orbit toroid.

Finally, assume that $\LL_{0} = s\bcln[n-1]$. As before, take $u \in \Pi^{\perp} \sm \{o\}$ closest to $o$, say $u=de_{n}$ for some $d \in \bN$. Since $\{2e_{1}, \dots , 2e_{n-2}, e_{1} + \cdots + e_{n-1}\}$ is a basis for $\LL_{0}$, by Lemma\nobreakspace \ref {lem:refLattices} there exist $\alpha_{1}, \dots, \alpha_{n-1}$ such that
\begin{align*}
	w &= \frac{s}{2} \Big( \alpha_{1} (2e_{1}) +  \dots + \alpha_{n-2}(2e_{n-2}) + \alpha_{n-1}(e_{1} + \cdots +e_{n-1}) \Big) + \frac{d}{2}e_{n}\\
	&=  \frac{s}{2} \Big( (2\alpha_{1}+\alpha_{n-1})e_{1} + \cdots + (2\alpha_{n-2}+\alpha_{n-1})e_{n-2}+   \alpha_{n-1} e_{n-1}\Big)+ \frac{d}{2}e_{n}
\end{align*}
satisfies that $\LL = \trasLat$. For $i \in \{1, \dots, n-1\}$ let $w_{i}$ the image of $w$ under the reflection on the hyperplane $x_{i} = 0$. If $\alpha_{n-1} = 1$ then $w - w_{n-1} = s e_{n-1}$, which contradicts that $\LL_{0} = s\bcln[n-1]$; then $\alpha_{n-1} = 0$. Let $i \in \{1, \dots, n-2\}$ such that $\alpha_{i} = 1$ and let $S_{i}$ the reflection $(x_{1}, \dots, x_{n}) \mapsto (x_{1}, \dots, x_{i-1}, x_{n-1}, x_{i+1}, \dots x_{n-2}, x_{i}, x_{n} )$, then $w-w_{i}S_{i} = s(e_{i} - e_{n-1})$ which for $n \geq 4$ does not belong to $s \bcln[n-1]$. This is a contradiction. Therefore if $\LL_{0} = s \bcln[n-1]$, then $\LL$ must be a vertical translation lattice.

The previous discussion completes the classification of cubic toroids with $n$ flag-orbits. We sumarize such discussion in the following result.

\begin{theorem}\label{thm:n-orbit}
	Let $n \geq 4$. Every equivelar $n$-orbit toroid of rank $(n+1)$ may be described as a toroid $\{4, 3^{n-2}, 4\}_{\bLL}$ where $\bLL$ is the lattice group associated to one of the following lattices:
	\begin{enumerate}
		\item $\bigcup_{k \in \bZ} \left( s \cln[n-2] + k(de_{n}) \right)$ for some $s, d \in \bN$, $s \neq d$.
		\item $\bigcup_{k \in \bZ} \left( s \fcln[n-3] + k(de_{n}) \right)$ for some $s, d \in \bN$.
		\item $\bigcup_{k \in \bZ} \left( s \bcln[n-1] + k(de_{n}) \right)$ for some $s, d \in \bN$.

		\item $\bigcup_{k \in \bZ} \left(  \cln[n-2] + k(\frac{s}{2} e_{1} +\frac{s}{2} e_{2}+  \cdots+ \frac{s}{2} e_{n-1} +  \frac{d}{2}e_{n} ) \right)$ for some $s, d \in \bN$, $s$ and $d$ even, and $d \neq s$.
		\item $\bigcup_{k \in \bZ} \left( s \fcln[n-3] + k(s e_{1}+ \frac
		{d}{2}e_{n})) \right)$ for some $s, d \in \bN$, $d$ even, and $d \neq 2s$.
	\end{enumerate}

\end{theorem}

%% file: nonCubic.tex
So far we have classified the cubic few-orbit  toroids arising from the $\{4,3^{n-2},4\}$ tessellation of the euclidean $n$-space. If $n\neq 2,4$ 
every equivelar toroid is a quotient of such tiling. The equivelar toroids for $n=2,3$ are fully described in \cite{hubardOrbanicPellicerWeissEquiv4Toroids}. Nevertheless, for $n=4$ the euclidean space admits a tessellation of type  $\{3,4,3,3\}$ by 24-cells, along with its dual, of type  $\{3,3,4,3\}$, by regular $4$-cross polytopes.

In this section we classify the equivelar toroids that are quotients of the regular tessellation of $\E[4]$ of type $\{3,3,4,3\}$. We will follow the same ideas we used in the classification of the cubic toroids. Recall that the vertex set of $\{3,3,4,3\}$ is the lattice $\bclFour$, and therefore its automorphism group has the form $\T(\cU)\rtimes \Go(\cU)=\langle R_0,\ldots,R_n\rangle$. Its vertex stabilizer $\langle R_1,\ldots,R_4\rangle$, with $R_1,R_2,R_3,R_4$ defined as in Table\nobreakspace \ref {tab:tessellations}, is  isomorphic to $[3,4,3]$, the automorphism group of the $24$-cell $\{3,4,3\}$. To proceed as in the case of the cubic toroids, we need some facts about the structure of $[3,4,3]$ that we present in the following paragraphs.

First, observe that $(\pm2,0,0,0),\ldots, (0,0,0,\pm2)$ and $(\pm 1,\pm 1,\pm 1,\pm 1)$ are the vertices of a 24-cell whose automorphism group is precisely the group $\langle R_1,R_2,R_3,R_4\rangle=[3,4,3]$.
Furthermore, they are partitioned into the sets $\mathcal{O}_0$, consisting of the points $(\pm 2,0,0,0)$, $(0,\pm 2,0,0)$, $(0,0,\pm 2,0)$ and $(0,0,0,\pm 2)$; $\mathcal{O}_1$, the set of points of the form $(\pm 1, \pm 1, \pm 1, \pm 1)$ with an even number of $-1$ entries, and $\cO_2$ the set of points of the form $(\pm 1, \pm 1, \pm 1, \pm 1)$ with an odd number of entries equal to $-1$. 

The vertex set of the 24-cell described above is preserved by $\stabcube$, implying that $\stabcube[4]\leq [3,4,3]$. Furthermore, the partition $\{\cO_{0}, \cO_{1}, \cO_{2}\}$ is preserved by every element of $[3,4,3]$. Also note that the elements of $[3,4,3]$ that fix every element of the partition are precisely those in $\abelnplus[4]\rtimes\symn[4]$. Hence $\abelnplus[4]\rtimes\symn[4]$ is a normal subgroup of $[3,4,3]$. 

Moreover, $\langle R_1,R_2\rangle\cap(\abelnplus[4]\rtimes\symn[4])=\{1\}$  and $|\abelnplus[4]\rtimes\symn[4]|=2^3\cdot 4!=\frac{|[3,4,3]|}{|\langle R_1,R_2\rangle|}$. Therefore, $[3,4,3]=(\abelnplus[4]\rtimes\symn[4])\rtimes\langle R_1,R_2\rangle$.

Finally, it is not hard to see that $\langle R_1,R_2\rangle$ acts as the full symmetric group on the elements of the partition, this observation will prove to be very useful to compute the symmetry types of the $3$-orbit non-cubic toroids.

%\subsection{The structure of $[3,4,3]$}
%From Table\nobreakspace \ref {tab:tessellations}, we get that $(0,0,0,0)$ and $(\pm 1,\pm 1,\pm 1,\pm 1)$, $(0^i,\pm 2,0^{3-i})$, for $0\leq i\leq 3$,  are vertices of $\{3,3,4,3\}$. Therefore, the convex hull of $(\pm\frac{1}{2},\pm\frac{1}{2},\pm\frac{1}{2},\pm\frac{1}{2})$, $(0^{i},\pm 1,0^{3-i})$, for $0\leq i\leq 3$, the standard cross polyope, is a vertex figure of this tessellation that we will refer to simply as $\{3,4,3\}$. This implies that $\stabcube[4]$ is a subgroup of $[3,4,3]$; furthermore, $\mathcal{O}:=\{(0^{i},\pm 1,0^{3-i})\,:\, 0\leq i\leq3\}$ is the vertex set of a four dimensional cross polytope and its orbits under $\langle R_1R_2\rangle$, stabilized by $\abelnplus[4]\rtimes\symn[4]$ are the  full vertex set of $\{3,4,3\}$. Observe that $\mathcal{O}(R_1R_2)$ and $\mathcal{O}(R_1R_2)^2$ and are the points of the form $(\pm\frac{1}{2},\pm\frac{1}{2},\pm\frac{1}{2},\pm\frac{1}{2})$ with an even and an odd number of minus signs, respectively.

%For $|\abelnplus[4]\rtimes\symn[4]|=\frac{|[3,4,3]|}{|\langle R_1,R_2\rangle|}$, from the previous paragraph, we obtain that $[3,4,3]=(\abelnplus[4]\rtimes\symn[4])\rtimes\langle R_1,R_2\rangle$.

A straightforward implementation of the {\tt ConjugacyClassesSubgroups} procedure in GAP System \cite{GAP4} can be used to find representatives of conjugacy classes of subgroups of index $2$, $3$ and $4$ that contain $\minid$. Those representatives are listed in Table\nobreakspace \ref {tab:lowindex}.

\begin{table}
\begin{center}
\begin{tabular}{|c|c|}
 \hline
   Index & Representatives \\ \hline
   \multirow{3}{*}{2} & $[3,4,3]^+$, the rotational subgroup of $[3,4,3]$,\\ \cline{2-2}
		     &$(\abelnplus[4]\rtimes\altn[4])\rtimes\langle R_1,R_2\rangle$,\\
		     \cline{2-2}
		     &$(\abelnplus[4]\rtimes\symn[4])\rtimes\langle R_1R_2\rangle$,\\ \hline
  \multirow{2}{*}{3} & $(\abelnplus[4]\rtimes\symn[4])\rtimes\langle R_2\rangle$ \\ \cline{2-2}
			& $(\abelnplus[4]\rtimes\died)\rtimes\langle R_1,R_2\rangle$\\ \hline
   4 & $(\abelnplus[4]\rtimes\altn[4])\rtimes\langle R_1R_2\rangle$ \\ \hline
  
 \end{tabular}
\caption{Representatives of conjugacy classes of small index in $[3,4,3]$}\label{tab:lowindex}
\end{center}
\end{table}

Now we can study the lattices contained in $\bclFour$ that are invariant under the representatives of the aforementioned conjugacy classes.

\subsection{\texorpdfstring{$2$}{2}-orbit non-cubic toroids} 

Let $G$ be one of the index 2 subgroups of $[3,4,3]$ listed above and $\LL \subseteq \bclFour$ be invariant under $G$. Since $(\abelnplus[4]\rtimes\altn[4])\rtimes\langle R_1R_2\rangle \leq G$, we have $\abelnplus[4]\rtimes\altn[4]\leq G$. Thus, by Lemma\nobreakspace \ref {lem:laticesC2plusAlt}, $\LL$ is of the form $s\LL_i$, $i\in\{0,1\}$, or $s\LL_{(1^k,0^{4-k})}$, for $k\in\{1,2,4\}$ and $s\in\mathbb{Z}^+$. 

Observe that the lattices  $s\LL_0$, $ s\LL_1$ and $s\bclFour$ are subsets of $\bclFour$ for every $s\in\mathbb{Z}^+$ , but $s\regLat$, $k\in\{1,2\}$ are contained in $\bclFour$ if and only if $s$ is even. Furthermore, for  every $s\in\mathbb{Z}^+$, $s\LL_{(1,0,0,0)}$ is not invariant under $(\abelnplus[4]\rtimes\altn[4])\rtimes\langle R_1R_2\rangle$, since $(s,0,0,0)R_1R_2=(\frac{s}{2},\frac{s}{2},\frac{s}{2},-\frac{s}{2})\notin s\LL_{(1,0,0,0)}$. Also, it is known that $2s\bclFour$ and $2s\LL_{(1,1,0,0)}$ are invariant under the full group, $[3,4,3]$, (see \cite[Section 6E]{ARP}). Therefore, $\{3,3,4,3\}_{\bLL}$ is regular for such lattices.

Let us consider first the group $\abelnplus[4]\rtimes\altn[4])\rtimes\langle R_1,R_2\rangle$. Since $\abeln[4] \rtimes \altn[4]$ is one of its subgroups, the proof of Theorem\nobreakspace \ref {thm:cubic2orbits} provides that the invariant lattices under this group must be of the form $s\regLat$, $k\in\{1,2,4\}$, but those lattices were ruled out in the previous paragraph. Therefore, there are no toroids $\{3,3,4,3\}_\bLL$ with two flag orbits and $\LL$ invariant under $\abelnplus[4]\rtimes\altn[4])\rtimes\langle R_1,R_2\rangle$.

Secondly, a toroid $\{3,3,4,3\}_\bLL$ with $\LL$ preserved by $[3,4,3]^+$ would be chiral, but such toroids were proven to be inexistent in \cite[Theorem 9.1]{mcMullenShulte_HigherToroidalRP}.

Finally, by the preceding discussion, if a toroid $\{3,3,4,3\}_\bLL$ is such that $\LL$ is preserved by $(\abelnplus[4]\rtimes\symn[4])\rtimes\langle R_1R_2\rangle$, then $\LL$ ought to be an integer multiple of the isometric lattices $\LL_0$ and $\LL_1$, known to be invariant under the action of $\abelnplus[4]\rtimes\symn[4]$. An easy computation shows that aforesaid lattices are also invariant under $R_1R_2$, from which it follows that for every positive integer $s$, $s\latOdd$ and $s\latEven$ are invariant under the action of $(\abelnplus[4]\rtimes\symn[4])\rtimes\langle R_1R_2\rangle$. Consequently the toroid $\{3,3,4,3\}_\mathbf{\LL}$, with $\bLL$ equal to $s \mathbf{\latOdd}$ or $s \mathbf{\latEven}$ for some integer $s$,  has two flag-orbits.

If $\UoverL$ is a 2-orbit toroid of type $\{3,3,4,3\}$, a similar analysis to the one made in Section\nobreakspace \ref {sec:cubic} shows that $ \UoverL $  is in class $2_{\{3,4\}}$, as $\LL_1$ is not preserved by $R_0$, $R_1$ or $R_2$. By duality, every $2$-orbit toroid of type $\{3,4,3,3\}$ is in class $2_{\{0,1\}}$.  We summarize the previous discussion in the following result.

\begin{theorem}
Let $\mathbf{\LL}$ be a lattice group of the form $\mathbb{s\latOdd}$, with $s$ an integer and $\latOdd$ the lattice group generated by the translations with respect to the vectors $(-1,1,1,1)$, $(1,-1,1,1)$, $(1,1,-1,1)$ and $(1,1,1,-1)$. Then the toroid $\{3,3,4,3\}/\mathbf{\LL}$ is a $2$-orbit toroid in class $2_{\{3,4\}}$ and the toroid$\{3,4,3,3\}/\mathbf{\LL}$ is a $2$-orbit toroid in class $2_{\{0,1\}}$. Furthermore, every non-cubic equivelar toroid with two flag orbits is of this form.
\end{theorem}

As every sublattice of $\bclFour$ that is invariant under $(\abelnplus[4]\rtimes\altn[4])\rtimes\langle R_1R_2\rangle$ is invariant under one of the subgroups that yield either to a regular or $2$-orbit toroid, we can state the following theorem.

\begin{theorem}
There are no $4$-orbit toroids of type $\{3,3,4,3\}$ or type $\{3,4,3,3\}$.
\end{theorem}

\subsection{\texorpdfstring{$3$}{3}-orbit non-cubic toroids}

For the index 3 subgroups, observe that $(\abelnplus[4]\rtimes\symn[4])\rtimes\langle R_2\rangle=\stabcube[4]$. Consequently, the only lattices invariant under $(\abelnplus[4]\rtimes\symn[4])\rtimes\langle R_2\rangle$ are of the form $s\LL_{(1^k,0^{4-k})}$, for $k\in\{1,2,4\}$. From the previous analysis we derive that $2s\LL_{(1,0,0,0)}$ is the only family of lattices preserved under $(\abelnplus[4]\rtimes\symn[4])\rtimes\langle R_2\rangle$ that is not preserved under $[3,4,3]$.

For the other representative of the index 3 subgroups, we get that every lattice invariant under $(\abelnplus[4]\rtimes\died)\rtimes\langle R_1,R_2\rangle$ should be preserved by $\abelnplus[4]\rtimes\died$, the only instances of those not preserved also by $[3,4,3]$ are of the form $2s(\LL_{(1,1)}\times\LL_{(1,1)})$. 

Thus, we have proved the following result.

\begin{theorem}\label{thm:noncubic}
 The non-cubic equivelar toroids with $3$ flag-orbits are exactly those of the form $\{3,3,4,3\}/\mathbf{\LL}$ or $\{3,4,3,3\}/\mathbf{\LL}$, with $\bLL$ either an even integer multiple of the lattice group $\mathbf{\LL_{(1,0,0,0)}}$ or an even integer multiple of the lattice group generated by the translations with respect to the vectors $(1,0,1,0)$, $(1,0,-1,0)$, $(0,1,0,1)$ and $(0,1,0,-1)$.
\end{theorem}

%Ahí me editas, pero creo que las ideas son esas.

%If $\LL\subseteq\bclFour$ is a lattice preserved under $[3,4,3]^+$, then $\{3,3,4,3\}/\LL$ must be a chiral toroid, but it is known that such toroids do not exist for ranks other than 3. Therefore, if a 2-orbit politope $\{3,3,4,3\}/\LL$ exists, then $\LL$ must be invariant under $\abelnplus[4]\rtimes\altn[4])\rtimes\langle R_1,R_2\rangle$ or $\abelnplus[4]\rtimes\symn[4])\rtimes\langle R_1R_2\rangle$.

%Observe that both possibilities are mutually exclusive, since in that case $\LL$ must be preserved under $\stabcube$, which implies that $\LL$ is of the form $\LL_{(1^k,0^{4-k})}$, for $k\in\{1,2,4\}$, but $\{3,3,4,3\}/\LL$ is regular\footnote{No estoy seguro de esto, la notación ce \cite{ARP} de $\LL_{cosas}$ es ligeramente distinta a la de nosotros cuando el grupo no es el del cubo. Se me hace que $\bclFour$ no es invariante respecto a todo el grupo.}. in all these cases \cite{ARP} 

%% file: symType.tex
In this section we determine the symmetry type of each family of few-orbit toroids classified in Sections\nobreakspace \ref {sec:cubic} and\nobreakspace  \ref {sec:nocubic}. Following \cite{cunninghamDelRioHubardToledo_stgPolytopesManiplexes}, the \emph{symmetry type graph} $\mathcal{T}(\cU/\bLL)$ of a toroid $\cU/\bLL$, is the labeled pre-graph (that is, semi-edges and and multiple edges are allowed) whose vertex-set is the set of orbits of flags of $\cU/\bLL$ and such that there is an edge (or semi-edge) labeled $i$ between orbits $\cO_{1}$ and $\cO_{2}$ if and only if there exists a flag $\Phi$ such that $\Phi \in \cO_{1}$ and $\Phi^{i} \in \cO_{2}$. Observe that since $(\Phi^{i})S = ( \Phi S )^{i}$ for every flag $\Phi$ and every automorphism $S$, the symmetry type graph does not depend on the representatives of the flag-orbits. We shall agree on using a semi-edge labeled $i$ instead of a loop whenever $\Phi$ and $\Phi^{i}$ belong to the same orbit. Note that this definition is slightly different of that of \cite{cunninghamDelRioHubardToledo_stgPolytopesManiplexes}, however it is easy to see that both are equivalent.

According to this definition, the symmetry type graph of a regular $n$-toroid consists of only one vertex and $n$ semi-edges labeled with $\{0, \dots, n-1\}$. If an $n$-toroid is in class $2_{I}$ for $I \subseteq \{0, \dots, n-1\}$, then its symmetry type graph is composed of only two vertices with edges whose labels lie in $\{0, \dots, n-1\} \sm I$ and semi-edges on each vertex labeled with the elements of $I$. In this sense, the symmetry type graph generalizes the notion of toroids in class $2_{I}$ for toroids with more than two flag-orbits. Observe that the symmetry type graph of the dual of a $n$-toroid its just the pre-graph with the same vertex set as the symmetry type graph of the toroid and with an edge labeled $n-i-1$ for every $i$-labeled edge.

Symmetry type graphs describe not only the number of flag-orbits of a toroid, but also the local arrangement of the orbits. In order to determine the symmetry graph of few-orbit toroids we use the following result, which is essentially  \cite[Proposition 1]{cunninghamDelRioHubardToledo_stgPolytopesManiplexes} on the language of toroids.

\begin{lemma}\label{lem:stgFaceTransitivity}
	Let $\cU/\bLL$ be a toroid with symmetry type graph $\mathcal{T}(\cU/\bLL)$. Let $\mathcal{T}^{i}(\cU/\bLL)$ be the the subgraph of $\mathcal{T}(\cU/\bLL)$ obtained by erasing the edges labeled $i$ of $\mathcal{T}(\cU/\bLL)$. Then $\cU/\bLL$ is $i$-face-transitive if and only if $\mathcal{T}^{i}(\cU/\bLL)$ is connected.
\end{lemma}

Recall that the group of automorphisms of a toroid $\UoverL$ is the group $\normGU(\bLL)/\bLL$. If $\cU$ is a regular tessellation of $\E$ then, up to duality, $\T(\cU)$ acts transitively on the vertices of  $\,\cU$. This implies that all the flag-orbits of $\,\UoverL$ occur on the base vertex of $\UoverL$. Furthermore, with the correspondence introduced in Lemma\nobreakspace \ref {lem:toroids} between $\normGU(\bLL)$ and a subgroup $N$ of the vertex stabilizer of $\G(\cU)$, the configuration of flag-orbits of $\UoverL$ around the base vertex is the same as the configuration of orbits of flags of $\cU$ containing the base vertex under the action of $N$. We use Lemma\nobreakspace \ref {lem:stgFaceTransitivity} and the previous observation, together with the results of \cite{cunninghamDelRioHubardToledo_stgPolytopesManiplexes} on symmetry type graphs of highly symmetric maniplexes to determine the symmetry type graph of the few-orbit toroids.

\subsection{Cubic toroids}

Symmetry type graphs of regular and $2$-orbit toroids were already described above. Since $k$-orbit $(n+1)$-toroids do not exist for $2<k<n$ unless $n=4$ we only need to determine the symmetry type graphs of $3$-orbit $(4+1)$-toroids and $n$-orbit $(n+1)$-toroids for $n$ even and $n \geq 4$.

Recall that if $\UoverL$ is a cubic toroid, then $\autUoverL = \normGU(\bLL)/ \bLL = (\T(\cU) \rtimes N) / \bLL $ for a certain group $N$ with $N \leq S$, where $S$ denotes the stabilizer of the vertex $o$ under $\G(\cU)$. Note that the symmetry type graph depends only on $N$, in particular all $(n+1)$-toroids with $n$ orbits share the same symmetry type graph.    

Let $\UoverL$ a $(4+1)$-toroid with $3$ flag-orbits. According to Lemma\nobreakspace \ref {lem:4cubicIndex3Groups}, the group $N$ is the group $\abeln[4] \rtimes \died$, where $\died$ acts on the coordinates of $\E[4]$ as the dihedral group. Note that since $\T(\cU) \leq \normGU(\bLL)$, $\autUoverL$ acts transitively on vertices and on cells of $\UoverL$. Vertex transitivity of $\UoverL$ implies that every edge of $\cU$ is in the same orbit as one joining the vertex $o$ with the vertex $\pm e_{i}$ for some $i \in \{1,2,3,4\}$. Since $\abeln[4] \leq N$ we may assume that the sign is positive. However, $\died$ acts transitively on the points $e_{i}$ for $i \in \{1,2,3,4\}$. This implies that $\UoverL$ is edge transitive. A similar argument can be used to prove that $\UoverL$ is transitive on rank-$3$ faces. Since there is no $3$-orbit $i$-face
transitive toroid for every $i\in \{0, \dots, 4\}$ (see \cite[Theorem 1]{cunninghamDelRioHubardToledo_stgPolytopesManiplexes}), then the symmetry type graph of $\UoverL$ must be that of Figure\nobreakspace \ref {fig:stg3orb1}.

To determine the symmetry type graph of an $n$-orbit $(n+1)$-toroid we shall change slightly the approach. Instead of looking at the action of the vertex-stabilizer of $\autUoverL$ we will use the stabilizer of a cell.  One way to imagine this is no to think $o$ as the vertex of $\UoverL$ but as the center of a cell. (alternatively we may just take the dual of $\UoverL$, which is isomorphic to $\UoverL$ since $\cU$ is self-dual). With this in mind we may think that the group $N = \abeln \rtimes \symn[n-1]$ acts on a cell.

Note that even though our results regarding cubic $n$-orbit $(n+1)$-toroids are for $n \geq 4$, the ideas apply as well when $n=3$. In particular these give a way to classify $4$-toroids on class $3$ listed in Tables 2 and 3 of \cite{hubardOrbanicPellicerWeissEquiv4Toroids}. Observe that the symmetry type graph of such toroids is that of Figure\nobreakspace \ref {fig:stgCubic2} for $n=3$ (see \cite[Fig. 4]{hubardOrbanicPellicerWeissEquiv4Toroids}). 

The observation in the previous paragraph allows us to give the following inductive argument to determine the symmetry type graph of any $(n+1)$-toroid with $n$ orbits. Let $n \geq 3$ and let $\UoverL$ be an $(n+1)$-toroid with $n$ orbits. We will show that the symmetry type graph of $\UoverL$ is that of Figure\nobreakspace \ref {fig:stgCubic2}. The case $n=3$ is explained above; assume then that $n \geq 4$. As said before, we shall think that $o$ is the center of a cell $C_{o}$ whose edges are segments of lenght $1$. Observe that the group $N = \abeln \rtimes \symn[n-1] $ acts transitively on the flags containing the facet $F$ of $C_{o}$ in the hyperplane $x_{n}=\frac{1}{2}$. Moreover, $N$ acts transitively on the set of facets of $C_{o}$ contained in the hyperplanes $x_{i}=\pm \frac{1}{2}$ for $i \neq n$. This implies that all the other flag-obits of $\UoverL$ have representatives on the flags containing the facet $F R_{n-1}$ of $C_{o}$. This facet is contained on the hyperplane $x_{n-1} = \frac{1}{2}$. However, the stabilizer of such face under the action of $N$ is precisely $\abeln[n-1] \rtimes \symn[n-2]$, where $\abeln[n-1]$ denotes the group generated by the reflections on all the coordinate hyperplanes $x_{i} =0$ for $i \neq n-1$, and $\symn[n-2]$ denotes the point-wise stabilizer of the $(n-1)$-th and $n$-th coordinates. By inductive hypothesis, the arrangement of flag-orbits of $F R_{n-1}$ induce the graph on Figure\nobreakspace \ref {fig:stgCubic2} with $n-1$ vertices, implying that the symmetry type graph of $\UoverL$ is that of Figure\nobreakspace \ref {fig:stgCubic2}. Those results are summarized on the following theorem.

\begin{theorem}\label{thm:stgCubic}
The symmetry type graph of a cubic $(4+1)$-toroid with three flag orbits is the graph shown in Figure\nobreakspace \ref {fig:stg3orb1}. If $n \geq 3$, the symmetry type graph of a cubic $(n+1)$-toroid with $n$ flag orbits is the graph with $n$ vertices shown in Figure\nobreakspace \ref {fig:stgCubic2}.  
\end{theorem}

\subsection{Non-cubic toroids}

Now we will determine the symmetry type graphs of the few-orbit toroids of type $\{3,3,4,3\}$ and $\{3,4,3,3\}$. As mentioned before, we can restrict our work to those toroids of type $\{3,3,4,3\}$; the symmetry type graph or the toroids of type $\{3,4,3,3\}$ may be obtained by duality, namely, interchanging the labels $0$ and $4$ as well as $1$ and $3$ of the symetry type graph of the toroids of type $\{3,3,4,3\}$.

 According to Theorem\nobreakspace \ref {thm:noncubic} there are two families of regular toroids of type $\{3,3,4,3\}$, one family of $2$-orbit toroids and two families of $3$-orbit toroids. The symmetry type graph of the regular toroids and toroids in class $2_{\{3,4\}}$ are already described above. It only remains to determine the symmetry type graphs of $3$-orbit toroids.
 
 Just as in the case of cubic toroids, we will determine the transitivity of the automorphism group on the set of $i$-faces of each toroid and use the results of \cite{cunninghamDelRioHubardToledo_stgPolytopesManiplexes} to determine its symmetry type. 
 
 Recall that the vertex set of the tessellation $\cU=\{3,3,4,3\}$ is a lattice. Since every translation of $\G(\cU)$ induces an automorphism of any toroid $\UoverL$, the group $\autUoverL$ acts transitively on vertices. 
 
 As before, since $\autUoverL = \normGU(\bLL)/\bLL = (\T(\cU) \rtimes K)/ \bLL $ for some group $K \leq \Go(\cU)$, to determine whether or not $\autUoverL$ acts transitively on the $i$-faces of $\UoverL$ for $i \in \{1,2,3,4\}$, it is enough to determine if $K$ acts transitively on the $i$-faces of $\cU$ that contain the vertex $o$. The $i$-faces of $\cU$ containing $o$ are in correspondence with the $(i-1)$-faces of the convex polytope $\{3,4,3\}$ with vertices $(2,0,0,0), \dots, (0,0,0,2)$ and $(\pm 1,\pm 1,\pm 1,\pm 1)$. We will use this correspondence in the following paragraphs.
   
 Following the notation of Table\nobreakspace \ref {tab:lowindex}, the group $K$ for the two families of $3$-orbit toroids are $(\abelnplus[4] \rtimes \symn[4]) \rtimes \langle R_{2} \rangle = \stabcube[4]$ and $(\abelnplus[4] \rtimes \died) \rtimes \langle R_{1}, R_{2} \rangle$.
 
\begin{figure}
\centering
\begin{subfigure}[b]{0.4\textwidth}
  \centering
	\includegraphics[width=\textwidth]{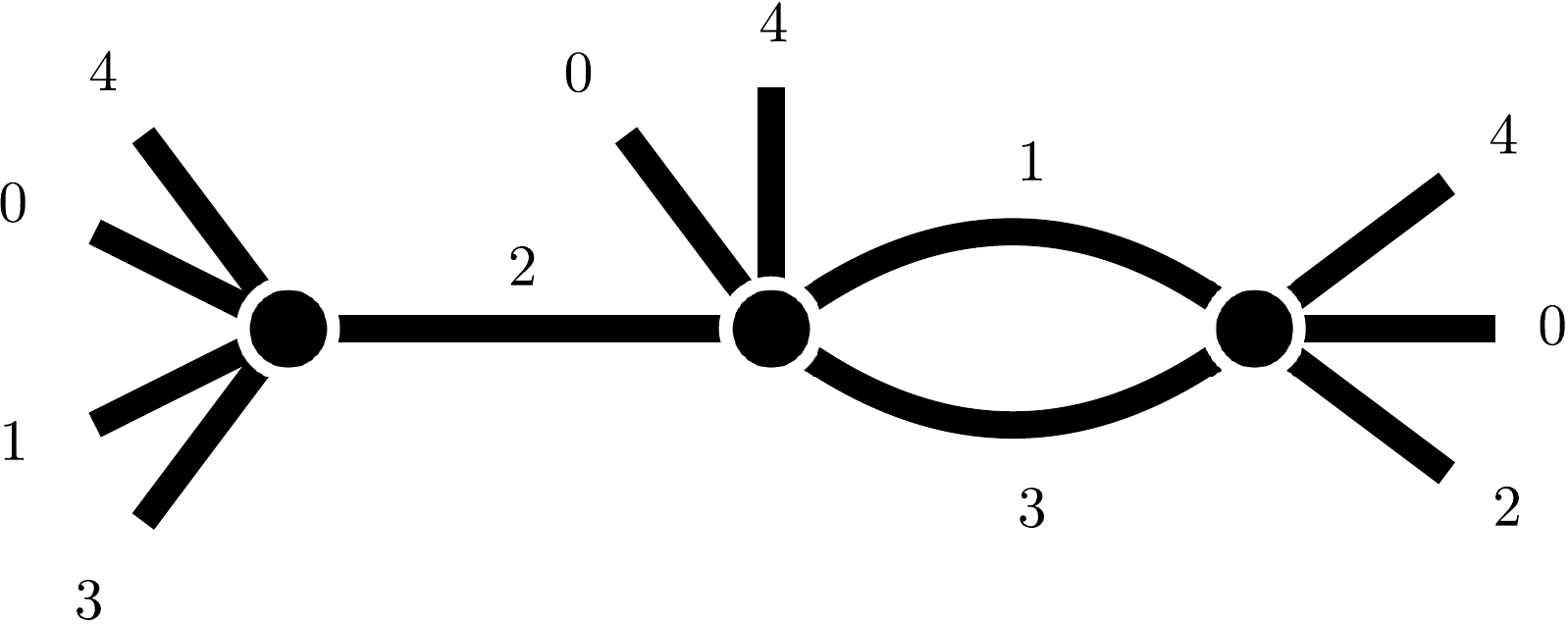}
  \caption{}\label{fig:stg3orb1}             
\end{subfigure}%

\begin{subfigure}[b]{0.4\textwidth}
  \centering
	\includegraphics[width=\textwidth]{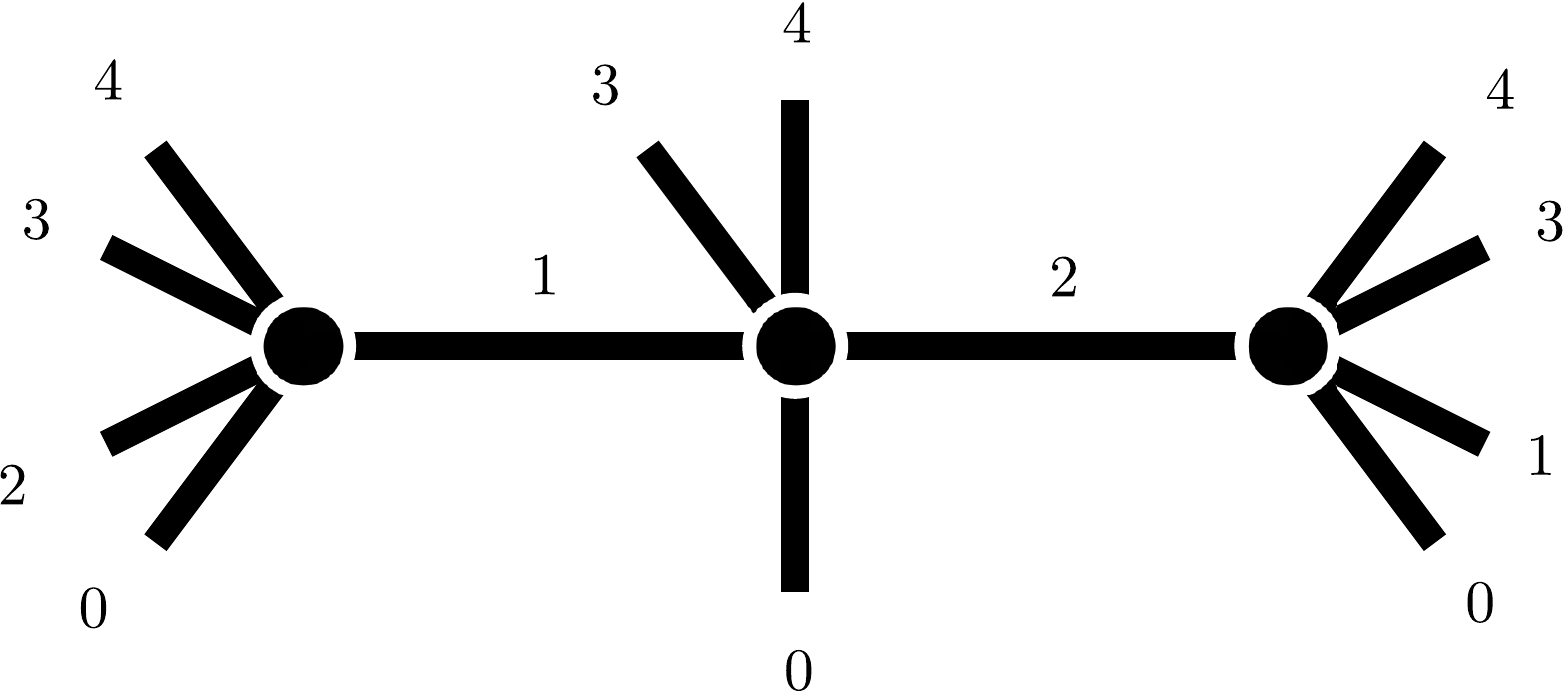}
  \caption{}\label{fig:stg3orb2}             
\end{subfigure}%

\caption{Symmetry type graphs of $3$-orbit $4$-toroids}
\end{figure}

\begin{figure}
 \centering
\begin{scriptsize}
	\includegraphics[width=\textwidth]{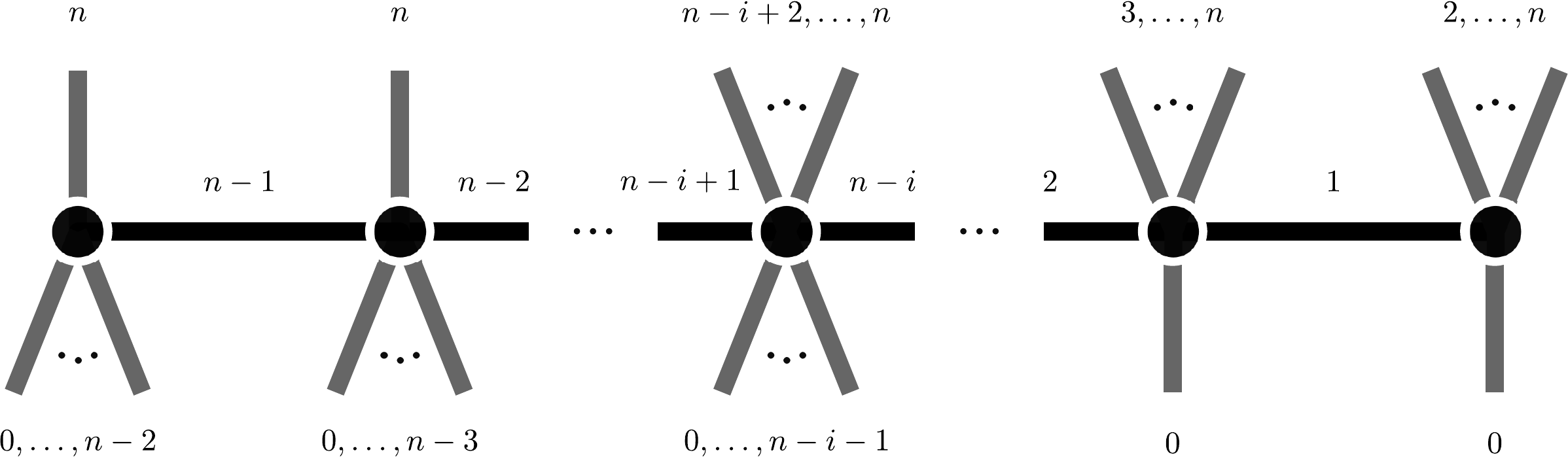}
  \caption{Symmetry type graph of cubic $n$-orbit $(n+1)$-toroids}\label{fig:stgCubic2} 
\end{scriptsize}

\end{figure}

 Let $\UoverL$ be a toroid of type $\{3,3,4,3\}$ such that the group $K$ described above is $\stabcube[4]$. The edges of $\cU$ containing $o$ are line segments $op$ where $p$ is one of the vertices of the $\{3,4,3\}$ mentioned above. Note that $K$ maps point of the form $(\pm 1, \pm 1, \pm 1, \pm 1)$ in points of the form $(\pm 1, \pm 1, \pm 1, \pm 1)$, implying that $\UoverL$ is not edge-transitive. The same argument proves that the triangle of $\cU$ determined by $o$, $(2,0,0,0)$ and $(1,1,1,1)$ and the triangle determined by $o$, $(1,1,1,1)$ and $(-1,1,1,1)$ belong to different orbits under the action of $K$, hence $\UoverL$ is not $2$-face-transitive. Therefore, the symmetry type graph of $\UoverL$ is that of the Figure\nobreakspace \ref {fig:stg3orb2} (see \cite[Proposition 3]{cunninghamDelRioHubardToledo_stgPolytopesManiplexes}).
 
 Note that the group $(\abelnplus[4] \rtimes \died) \rtimes \langle R_{1}, R_{2} \rangle$ acts transitively on the vertices of $\{3,4,3\}$, namely, the set $\{(2,0,0,0), \dots , (0,0,0,2)\} \cup \{(\pm 1,\pm 1,\pm 1,\pm 1)\}$. This implies that if $\UoverL$ is a toroid such that the group $K$ described above is $(\abelnplus[4] \rtimes \died) \rtimes \langle R_{1}, R_{2} \rangle$ , then $\UoverL$ is edge transitive. 
 
The stabilizer of $(2,0,0,0)$ under the action of $K$ has size $\frac{|K|}{24} = 16$ and contains the subgroup $H=(\langle E_{2}, E_{3}, E_{4}\rangle^{+} \rtimes \langle S_{(2,4)}\rangle ) \rtimes \langle R_{2} \rangle = \langle E_{2}, E_{3}, E_{4}\rangle \rtimes \langle S_{(2,4)}\rangle$, where $E_{i}$ denotes the reflection in the hyperplane $x_{i} =0$ and $S_{(2,4)}$ denotes the isometry that interchanges the second and the fourth coordinates. Note that $|H|=16$ and thus the stabilizer of $(2,0,0,0)$ is precisely $H$.

Assume that there is an isometry $S$ of $K$ mapping the triangle of $\cU$ determined by the vertices $o$, $(2,0,0,0)$ and $(1,1,1,1)$ to the triangle $o$, $(2,0,0,0)$, $(-1,-1,1,1)$. Note that $S$ cannot fix $(2,0,0,0)$ since there is no element of $H$ that maps $(1,1,1,1)$ to $(-1,-1,1,1)$. This implies that $S$ maps $(2,0,0,0)$ to $(1,1,1,1)$. Recall that the group $\langle R_{1}, R_{2} \rangle$ permutes the blocks of vertices $\cO_0=\{2e_{i} :   1 \leq i \leq 4 \}$, $\cO_1= \{(x_{1},x_{2}, x_{3}, x_{4})  :   |x_{i}| = 1, x_{1}x_{2}x_{3}x_{4} = 1 \}$ and $\cO_2 = \{(x_{1},x_{2}, x_{3}, x_{4})  :   |x_{i}| = 1, x_{1}x_{2}x_{3}x_{4} = -1 \}$. Note that $S$ maps a vertex of $\cO_{0}$ to a vertex of $O_{1}$, this implies that $S = T R_{1}$ with $T \in \abelnplus[4] \rtimes \died$. Observe that \[(1,1,1,1) = (2,0,0,0)S^{-1} = (2,0,0,0)R_{1}T^{-1} = (1,1,1,1)T^{-1},\] which implies that $T \in \died$. On the other hand \[(2,0,0,0)=(-1,-1,1,1)S^{-1} = (-1,-1,1,1)R_{1} T^{-1} = (0,-2,0,0)T^{-1}\] which contradicts that $T \in \died$. This implies that $K$ does not act transitively on the edges of $\{3,4,3\}$ hence $\UoverL$ is not transitive in $2$-faces

To see that $\UoverL$ is $3$-face-transitive it is enough to show that $K$ acts transitively on the triangles of $\{3,4,3\}$. We will show that $K$ acts transitively on the triangles of $\{3,4,3\}$. incident to $(2,0,0,0)$. Such triangles are in correspondence with the edges of the $3$-dimensional cube $Q$ determined by the vertices $(1, \pm 1, \pm 1, \pm 1)$. Let $H$ be as above, the stabilizer of $(2,0,0,0)$. The elements of $\langle E_{1}, E_{2}, E_{3} \rangle \leq H$ testify that the edges of $Q$ pointing in the same direction belong to the same orbit. The isometry $S_{(2,4)} \in H$ that permutes the second and the fourth coordinates, maps the edges in direction of $(0,2,0,0)$ to edges in direction $(0,0,0,2)$. Finally, the transformation $S: (x_{1}, x_{2},x_{3},x_{4}) \mapsto (x_{4}, x_{1},x_{2},x_{3})$ is an element of $\died$ that maps the edge of $Q$ determined by $(1,1,1,1)$ and $(1,-1,1,1)$ to the edge of $Q$ determined by $(1,1,1,1)$ and $(1,1,-1,1)$. This implies that all the edges of $Q$ belong to the same orbit under the action of $K$. Therefore, $K$ acts transitively in the triangles of $\cU$ incident to $o$ hence $\UoverL$ is $3$-face-transitive. 

The discussion in the previous paragraphs and Propositions 3 and 4 of \cite{cunninghamDelRioHubardToledo_stgPolytopesManiplexes} imply that the symetry type graph of a toroid $\UoverL$ whose automorphism group is $(\T(\cU) \rtimes K)/ \bLL$ with $K=(\abelnplus[4] \rtimes \died) \rtimes \langle R_{1}, R_{2} \rangle$ is the graph in Figure\nobreakspace \ref {fig:stg3orb1}. 

%% file: equivFewOrbitToroids.bbl
\begin{thebibliography}{10}

\bibitem{brehmKuhnel_EquivelarMapsTorus}
Ulrich Brehm and Wolfgang K\"uhnel.
\newblock Equivelar maps on the torus.
\newblock {\em European J. Combin.}, 29(8):1843--1861, 2008.

\bibitem{conderDobcsanyi_RegMapsSmallGenus}
Marston Conder and Peter Dobcs\'anyi.
\newblock Determination of all regular maps of small genus.
\newblock {\em J. Combin. Theory Ser. B}, 81(2):224--242, 2001.

\bibitem{conderRegMapsChar1to200}
Marston D.~E. Conder.
\newblock Regular maps and hypermaps of {E}uler characteristic {$-1$} to
  {$-200$}.
\newblock {\em J. Combin. Theory Ser. B}, 99(2):455--459, 2009.

\bibitem{coxeterRegularPolytopes}
H.~S.~M. Coxeter.
\newblock {\em Regular polytopes}.
\newblock Dover Publications, Inc., New York, third edition, 1973.

\bibitem{coxeterMoserGenandRelforDG}
H.~S.~M. Coxeter and W.~O.~J. Moser.
\newblock {\em Generators and relations for discrete groups}.
\newblock Springer-Verlag, New York-Heidelberg, third edition, 1972.
\newblock Ergebnisse der Mathematik und ihrer Grenzgebiete, Band 14.

\bibitem{cunninghamDelRioHubardToledo_stgPolytopesManiplexes}
Gabe Cunningham, Mar\'ia Del R\'io-Francos, Isabel Hubard, and Micael Toledo.
\newblock Symmetry type graphs of polytopes and maniplexes.
\newblock {\em Ann. Comb.}, 19(2):243--268, 2015.

\bibitem{dixonPermutationGroups}
John~D. Dixon and Brian Mortimer.
\newblock {\em Permutation groups}, volume 163 of {\em Graduate Texts in
  Mathematics}.
\newblock Springer-Verlag, New York, 1996.

\bibitem{GAP4}
The GAP~Group.
\newblock {\em GAP -- Groups, Algorithms, and Programming, Version 4.8.6},
  2016.

\bibitem{harleyMcMullenSchulte_SymTessESF}
Michael~I. Hartley, Peter McMullen, and Egon Schulte.
\newblock Symmetric tessellations on {E}uclidean space-forms.
\newblock {\em Canad. J. Math.}, 51(6):1230--1239, 1999.
\newblock Dedicated to H. S. M. Coxeter on the occasion of his 90th birthday.

\bibitem{hubardOrbanicPellicerWeissEquiv4Toroids}
Isabel Hubard, Alen Orbani{\'c}, Daniel Pellicer, and Asia~Ivi{\'c} Weiss.
\newblock Symmetries of equivelar 4-toroids.
\newblock {\em Discrete Comput. Geom.}, 48(4):1110--1136, 2012.

\bibitem{matousek_DiscreteGeometry}
Ji{\v{r}}{\'{\i}} Matou{\v{s}}ek.
\newblock {\em Lectures on discrete geometry}, volume 212 of {\em Graduate
  Texts in Mathematics}.
\newblock Springer-Verlag, New York, 2002.

\bibitem{mcMullenShulte_HigherToroidalRP}
Peter McMullen and Egon Schulte.
\newblock Higher toroidal regular polytopes.
\newblock {\em Adv. Math.}, 117(1):17--51, 1996.

\bibitem{ARP}
Peter McMullen and Egon Schulte.
\newblock {\em Abstract regular polytopes}, volume~92 of {\em Encyclopedia of
  Mathematics and its Applications}.
\newblock Cambridge University Press, Cambridge, 2002.

\bibitem{orbanicPellicerWeiss_MapsOperationsKorbitMaps}
Alen Orbani{\'c}, Daniel Pellicer, and Asia~Ivi{\'c} Weiss.
\newblock Map operations and {$k$}-orbit maps.
\newblock {\em J. Combin. Theory Ser. A}, 117(4):411--429, 2010.

\bibitem{ratcliffeHypManifolds}
John~G. Ratcliffe.
\newblock {\em Foundations of hyperbolic manifolds}, volume 149 of {\em
  Graduate Texts in Mathematics}.
\newblock Springer, New York, second edition, 2006.

\end{thebibliography}
